\documentclass{article}

\usepackage{natbib}
 \bibpunct[, ]{(}{)}{,}{a}{}{,}%
 \def\bibfont{\small}%
\usepackage[utf8]{inputenc}
\usepackage{lineno,hyperref}
\usepackage{amsmath, amsfonts, amssymb,float} 
\usepackage{mathpazo}
\usepackage{fullpage}
\usepackage{tabularx}
\usepackage{multirow}
\usepackage{graphicx}
\usepackage{subcaption}
\usepackage{nicefrac}
\usepackage{tabularx}
\usepackage{natbib}
\usepackage[ruled,vlined]{algorithm2e}
\usepackage{mathtools}
\DeclarePairedDelimiter{\ceil}{\lceil}{\rceil}

\providecommand{\keywords}[1]
{
  \small
  \textbf{\textit{Keywords---}} #1
}

\usepackage{authblk}
\title{Dynamic Pooled Capacity Deployment \\ for Urban Parcel Logistics}
\date{}
\author[1]{Louis Faug\`ere}
\author[2]{Walid Klibi}
\author[3,4]{Chelsea White III}
\author[1,3,5]{Benoit Montreuil}
\affil[1]{Physical Internet Center, School of Industrial \& Systems Engineering, Georgia Institute of Technology, Atlanta GA, USA}
\affil[2]{The Center of Excellence in Supply Chain (CESIT), KEDGE Business School, Bordeaux, France}
\affil[3]{Supply Chain \& Logistics Institute, Georgia Institute of Technology, Atlanta GA, USA}
\affil[4]{Schneider National Chair in Transportation and Logistics}
\affil[5]{Coca-Cola Material Handling \& Distribution Chair}
\begin{document}

\maketitle

\begin{abstract}
Last-mile logistics is regarded as an essential yet highly expensive component of parcel logistics. In dense urban environments, this is partially caused by inherent inefficiencies due to traffic congestion and the disparity and accessibility of customer locations. In parcel logistics, access hubs are facilities supporting relay-based last-mile activities by offering temporary storage locations enabling the decoupling of last-mile activities from the rest of the urban distribution chain. This paper focuses on a novel tactical problem: the geographically dynamic deployment of pooled relocatable storage capacity modules in an urban parcel network operating under space-time uncertainty. In particular, it proposes a two-stage stochastic optimization model for the access hub dynamic pooled capacity deployment problem with synchronization of underlying operations through travel time estimates, and a solution approach based on a rolling horizon algorithm with lookahead and a benders decomposition able to solve large scale instances of a real-sized megacity. Numerical results, inspired by the case of a large parcel express carrier, are provided to evaluate the computational performance of the proposed approach and suggest up to $28\%$ last-mile cost savings and $26\%$ capacity savings compared to a static capacity deployment strategy. \\
\keywords{Parcel Logistics, Urban Networks, Dynamic Deployment, Capacity Relocation, Capacity Pooling, Stochastic Optimization, Physical Internet}
\end{abstract}

\section{Introduction} \label{introduction}

\begin{figure}[H]
\centering
\fbox{\includegraphics[scale=0.30]{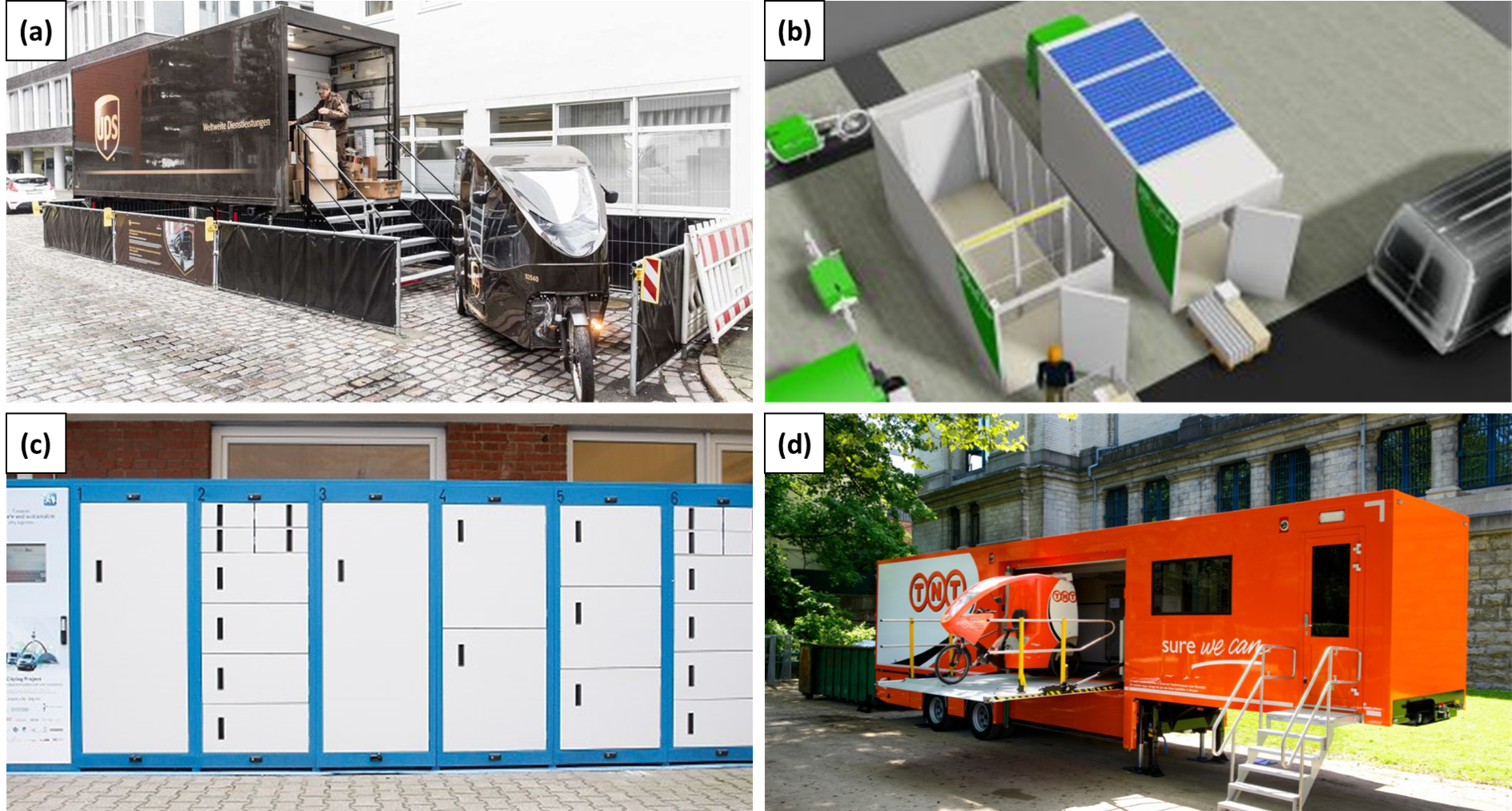}}
\caption{Examples of Access Hubs: (a) Trailer Micro-hub in Berlin by \cite{UPSMicroHub}, (b) Storage Shed by \cite{HubCompany}, (c) Transshipment Locker Bank by \cite{CityLogBentoBox}, and (d) Trailer by \cite{TNT:MobileAH}}
\label{fig: AHExample}
\end{figure}

Global urbanization, growth of e-commerce and the ever increasing desire for speed put pressure on the need for innovation in designing, managing and operating urban logistics systems in a sustainable and cost-efficient way. In 2018, $55\%$ of the world's population lived in urban areas (up to $82\%$ in North America). The \cite{united20182018} predict that global urbanization will reach $68\%$ by 2050, with an increasing number of megacities (cities of 10+M inhabitants). Increasing population density is a challenge for city logistics in terms of traffic congestion, vehicle type restrictions, limited parking spaces, expensive and rare logistic facility locations, and is further complex in megacities due to their extremely high density \citep{fransoo2017reaching}. For urban parcel logistics systems, the growth of e-commerce is currently one of the main challenges to tackle with an annual growth over $20\%$ on the 2017-2019 period, projected to be over $15\%$ until 2023 \citep{statista_2019_growth}. Online-retailing with goods being transported to consumers' homes increase the number of freight movements within cities while reducing the size of each shipment \citep{savelsbergh201650th} which makes first and last mile logistic activities harder to plan. Moreover, consumers' desire for speed (i.e. same-day delivery and faster) has yet to be met by online retailers \citep{statista_2019_speed}. With promises as fast as 1-hour delivery (e.g. Amazon Prime in select U.S cities), the cost of last-mile logistics becomes an ever more critical part of urban parcel logistics. These trends have been accelerated due to attempts to mitigate the impacts of the COVID-19 pandemic (e.g., sequestering in place), requiring companies to increase their last-mile delivery capabilities and to deal with the dramatic shift to online channels \citep{mitsloancovid19}.
\\
To tackle these challenges, a number of innovations have emerged from academia and industry. \cite{savelsbergh201650th} provide an overall view of recent innovations and modeling of solutions such as multi-echelon networks, dynamic delivery systems, pickup and delivery point networks, omni-channel logistics, crowd-sourced transportation and the integration of public and freight transportation networks. Many of these innovations are considered in the Physical Internet initiative, introduced in \citet{montreuil2011toward}, which seeks global logistics efficiency and sustainability by transforming the way physical objects are handled, moved, and stored by applying concepts from internet data transfer to real-world shipping processes. A conceptual framework on the application of Physical Internet concepts to city logistics was recently proposed in \citet{crainic2016physical}, in particular the concepts of pooling and hyperconnectivity in urban multi-echelon networks. As underlined by \citet{savelsbergh201650th}, city logistics problems integrating real-life features such as highly dynamic and volatile decision making environments, sharing principles or multi-echelon networks, offer a fertile soil for groundbreaking research.
\\
Inspired by the case of a large parcel logistics company operating in megacities, this paper examines a novel tactical optimization problem in urban parcel logistics. It consists in the dynamic deployment and relocation of pooled storage capacity in an urban parcel network operating under space-time uncertainty. It builds on the recent proposal of a hyperconnected urban logistics network structure  \citep{montreuil2018IMHRC} in line with the new challenges of the parcel logistics industry. The proposed network structure is based on the pixelization of urban agglomerations in unit zones (clusters of customer locations), local cells (cluster of unit zones) and urban areas (cluster of local cells). It is composed of three tiers of interconnected logistics hubs: gateway hubs (GH), local hubs (LH) and access hubs (AH) respectively designed to efficiently handle inter urban areas, inter local cells, and inter unit zones parcel flows. Beyond the realm of an urban agglomeration, the network of gateway hubs connects to a network of regional hubs (RH) covering entire blocks of the world (e.g. North America), and these regional hubs connect to a worldwide network of global hubs.
This paper focuses on access hubs which are small logistics hubs located at the neighborhood level within minutes of customers, enabling parcel transfer between different vehicle types temporarily holding parcels close to pickup and delivery points. Access hubs are to be used by logistics carriers, and not by consumers as smart lockers are. Access hubs can materialize in many forms including a parked trailer, a smart locker bank, or a storage shed as illustrated in Figure \ref{fig: AHExample}. Trailer based solutions like Figure \ref{fig: AHExample} (a) and (d) offer all-or-nothing mobile solutions, while capacity module based solutions like Figure \ref{fig: AHExample} (b) and (c) offer flexible capacity adjustment over time. The scope of this paper is a capacity module based solution.

Parcel logistics networks have undergone significant changes in the last 20 years, notably in urban contexts as seen in \cite{janjevic2020characterizing}, and have received an increasing attention in the academic literature. Strategic and tactical network design problems such as the ones examined by \cite{SmilowitzDaganzo,winkenbach2016enabling} approximate operations costs when designing and planning for multi-echelon networks. While network design problems are complex due to intricate interdependencies between strategic, tactical and operational decisions, continuum approximations (see \cite{ansari2018advancements}) are useful to capture operations complexity and take informed decisions. However, such approximations are typically used to estimate travel distance and cost, but not travel time and operations synchronization.
This paper considers access hubs to be modular in storage capacity similar to designs proposed in \cite{FaugereCIE}, such that capacity modules can be removed/added to adapt access hub's storage capacity. At the tactical level, capacity modules are to be deployed over a network of access hub locations; at the operational level, capacity modules are to be allocated to serve their access hub's need or neighboring locations via capacity pooling. In a dynamic setting, the associated problem can be related to a multi-period location-allocation problem which belongs to the NP-Hard complexity class \citep{manzini2008optimization}.
Once the capacity of the network of access hubs is adjusted, each access hub plays the role of a transshipment location between couriers performing pickup and delivery services within minutes of the access hub and riders transporting parcels between local hubs and a set of access hubs. Such transshipments require tight synchronization of the two tiers so as to provide efficient and timely pickup and delivery operations. This operational context mimics, on a hourly basis, a two-echelon pickup and delivery problem with synchronisation, which is a complex routing problem (see for instance \citet{Cuda2015}). Thus, the integration of operations in the tactical decision model leads to better capacity deployment decisions \citep{Klibi2016}, yet induces solvability challenges due to its combinatorial and stochastic-dynamic structure.

This paper studies a novel tactical optimization problem: the dynamic deployment of pooled storage capacity in an urban parcel network operating under space-time uncertainty. Its contribution is threefold: (1) the characterization of a new tactical problem for capacity deployment, motivated by dynamic aspects of urban parcel logistics needs, (2) the modeling of the access hub dynamic pooled capacity deployment problem as a two-stage stochastic program with synchronization of underlying operations through travel time estimates, and (3) the design of a solution approach based on a rolling horizon algorithm with lookahead and a benders decomposition able to solve large scale instances of a real-sized megacity. Numerical results, inspired by the case of a large parcel express carrier, are provided to evaluate the computational performance of the proposed approach and suggest up to $28\%$ last-mile cost savings and $26\%$ capacity savings compared to a static capacity deployment strategy.

Section \ref{Literature} summarizes the literature relevant to this type of problem, section \ref{Probdescription} describes the problem and proposes a mathematical modeling, section \ref{SolutionApproach} presents the proposed solution approach, section \ref{Results} provides an experiment setup and discusses results, and section \ref{Conclusion} highlights key takeaways and managerial insights, and identifies promising research avenues.

\section{Literature Review} \label{Literature}

Multi-echelon network for urban distribution have received a lot of attention in the academic literature (e.g. \cite{BenjellounTrends}, \cite{mancini2013multi}, \cite{janjevic2019integrating}), commonly using urban consolidation centers (UCC) to bundle goods outside the boundaries of urban areas. As reported in \cite{janjevic_development_2014}, several micro-consolidation initiatives have been proposed to downscale the consolidation effort by bundling goods at the neighborhood level using capillary networks of hubs located much closer to pickup and delivery points, defined as access hubs in the conceptual framework proposed by \cite{montreuil2018IMHRC}. Examples of such initiatives are satellite platforms (e.g. \cite{BenjellounTrends}), micro-consolidation centers (e.g. \cite{leonardi2012before}), mobile depots (e.g. \cite{marujo2018assessing}), and micro-depots \cite{stodick2019sustainable}.
Most of the focus has been on location and vehicle routing aspects (e.g. \cite{anderluh2017synchronizing} and \cite{enthoven2020two}) and cost and negative externalities assessment (e.g. \cite{verlinde2014does}, \cite{arvidsson2017ex}, \cite{marujo2018assessing}) in solutions using depots and cargo-bikes. To the best of the authors' knowledge, the dynamic management of access hub capacity for urban parcel logistics has not yet been studied in the academic literature. \\
The problem studied in this paper involves modular capacity relocation and a capacity pooling recourse mechanism impacting the operations of a two-echelon synchronization problem. In this section, a literature review on dynamic capacitated facility location problems and integrated urban network design problems is presented. \\
Dynamic facility location problems where systems are subject to varying environments (e.g. non-stationary demand) allow the relocation of facilities over time. \cite{arabani2012facility} provide a literature review on facility location dynamics, including problems with and without hub relocation. Innovations in the manufacturing industry have motivated the study of modular and mobile production and storage. \cite{marcotte2016introducing} have presented various threads of innovations such as distributed production, on-demand production, additive production, and mobile production, that would motivate and benefit from hyperconnected mobile production systems. \cite{marcotte2015modeling} and \cite{malladiererawhite} proposed mathematical modeling for production and inventory capacity relocation and allocation to manage multi-facility network facing stochastic demand. However, they examine small to medium networks far from the scale of urban parcel logistics networks and do not study operations synchronization. \cite{aghezzaf2005capacity} studied storage capacity expansion planning coupled to dynamic inventory relocation in the context of warehouse location allocation problems, but did not consider capacity reduction or relocation. \cite{ghiani2002capacitated}, \cite{melo2006dynamic}, and \cite{jena2015dynamic} modeled dynamic facility location problem where not only sites could be permanently or temporarily opened or closed, but also resized by adding or removing modular capacity. \cite{melo2006dynamic} proposed models capturing modular capacity shifts from existing to new facilities. However in these problems, capacity relocation is generally not managed jointly with capacity allocation or its impact on underlying operations. Dynamic facility location literature partially covers the tactical capacity relocation problem studied in this paper, but does not integrate underlying operations dynamics at the urban logistics scale. \\
Integrated network design problems typically deal with a combination of strategic decisions such as facility location, tactical decisions such as resource allocation and scheduling, and operational decision such as vehicle routing. The integration of these different levels of decisions can be found in two main problem classes: service network design problems and location routing problems. Service network design problems deal with the selection and scheduling of services such as hub operations, shipping lines and routing of freight (e.g. \cite{crainic2016service,hewitt2019scheduled}) while location routing problems combine facility location-allocation decisions with associated freight routing decisions. \cite{drexl2015survey} provide a recent survey of variants and extensions of the location routing problem. The dynamic location routing problem (\cite{francis2008period}) considering the assignment of demand to locations over multiple periods, is similar to the problem studied in this paper: it aims at minimizing network and routing costs over a multiperiod location and routing decision vector. However, multi-echelon location routing problems (e.g. \cite{crainic2004advanced,perboli2011two}) have only recently gathered attention in the literature. Although multi-echelon networks are relevant to postal and parcel delivery distribution systems (\cite{gonzalez2009n}) where fine time constraints and synchronization have become an essential consideration, most papers studying multi-echelon networks are concerned with the two-echelon case and ignore temporal aspects (\cite{drexl2015survey}). \\
When allowing inter-location capacity pooling, underlying operations described in section \ref{introduction} are impacted. Couriers perform pickup and delivery tours starting and ending in their reference access hub, while riders visit access hubs starting and ending their routes in their reference local hub. The impact of capacity pooling can be measured by modeling its impact on the route of parcels, couriers and riders. However, when taking decisions at the tactical level, explicitly modeling routes is not necessary. TSP and VRP continuous approximations have been introduced by \cite{Daganzo94, DaganzoBook} to embed operations in strategic and tactical logistics problems (e.g. \cite{EreraThesis}, \cite{franceschetti2017strategic}). A recent literature on variants of this approach can be found in \cite{ansari2018advancements}. \cite{SmilowitzDaganzo, winkenbach2016enabling,bergmann2020integrating} adapted these continuous approximations to the context of parcel express logistics to approximate distance traveled and cost. However, the aspect of synchronization using travel time continuous approximations has not yet been studied.
To the best of the authors' knowledge, this paper is the first to study a capacity relocation problem with the synchronization of two-echelon routing operations through travel time estimates.


\section{Problem Description and Formulation} \label{Probdescription}
\subsection{Business Context}\label{biz}
\begin{figure}[H]
\centering
\fbox{\includegraphics[scale=0.50]{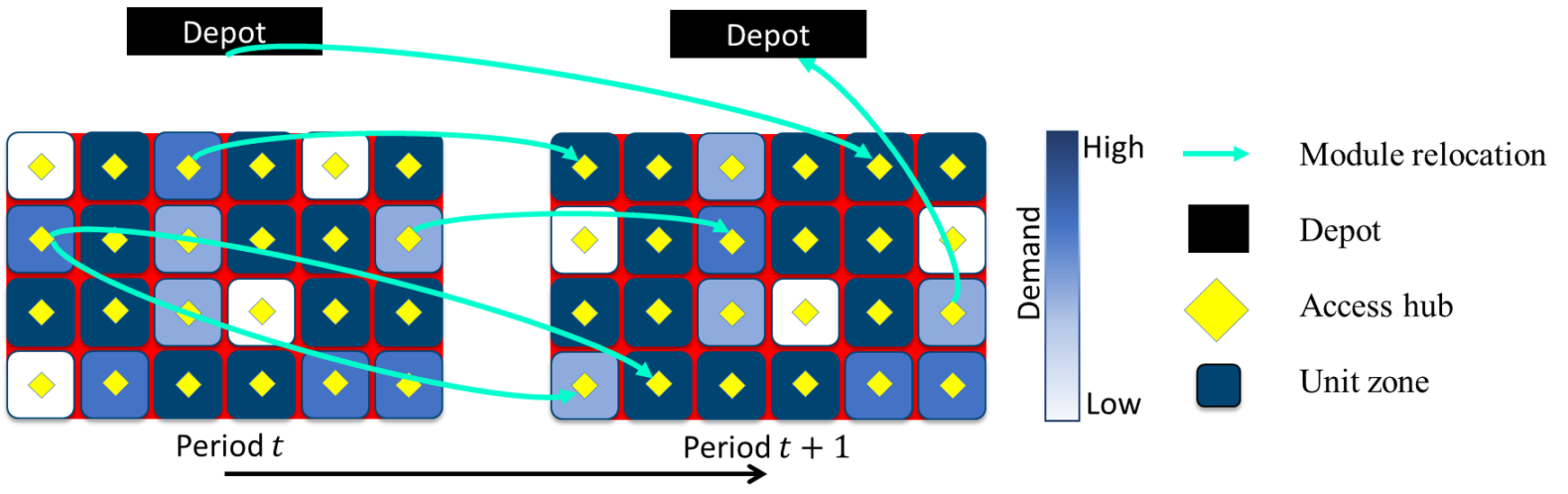}}
\caption{Illustration of Access Hub Capacity Relocation }
\label{fig: Relocation}
\end{figure}

A parcel logistics company provides pickup and delivery services to customers in a region covered by a network of access hubs. The network of access hubs may be dedicated to the parcel logistics provider, or shared between several companies as suggested by the concept of open networks in the Physical Internet. Figure \ref{fig: Relocation} provides a conceptual illustration of the network of access hubs and the relocation of capacity modules over two deployment periods. Once the network capacity is set, pickups from customers are dropped off by couriers in access hubs and will occupy a certain storage volume for some time until a rider  picks them up to perform outbound activities. To-be-delivered parcels are dropped off by riders in access hubs and will occupy a certain storage volume for some time until a courier picks them up to perform the delivery to customers.
To provide good service, the company must ensure that parcels flow rapidly and seamlessly between couriers and riders, which requires the sound management of storage capacity deployed in access hubs. Storage volume requirements vary depending on the fluctuation of demand for pickup and delivery services over time and are observed over a discrete set of operational periods (e.g. hourly). Access hubs are composed of modular storage units that can be assembled and disassembled relatively easily, enabling rapid relocation of storage capacity in the network. During each deployment period (e.g. week or day), storage capacity can be relocated within the network of access hubs, or to/from a depot where additional capacity modules are stored when not in use.
Figure \ref{fig: Relocation} illustrates demand variability and the relocation of capacity modules within the network of access hubs over two deployment periods. For instance, unit zones with increasing demand (and therefore increasing capacity requirements) from period $t$ to $t+1$ receive capacity module(s) from the depot of from locations that have decreasing capacity requirements (e.g. lower left unit zone in Figure \ref{fig: Relocation}).

\begin{figure}[H]
\centering
\fbox{\includegraphics[scale=0.50]{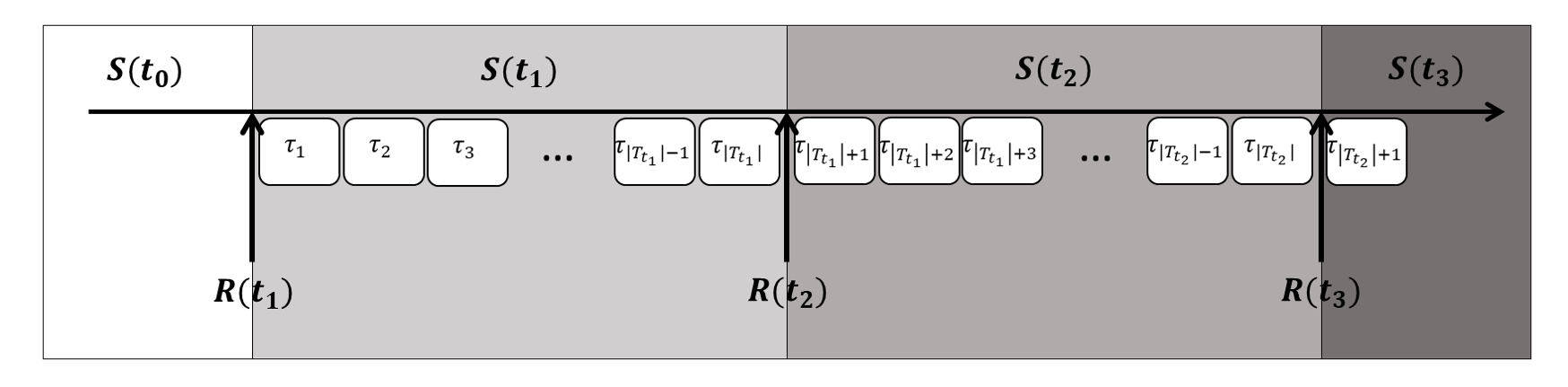}}
\caption{Timeline of the Hyperconnected Access Hub Network Dynamic Storage Capacity Deployment Problem}
\label{fig: Timeline}
\end{figure}

The relocation of capacity modules over the network adjusts the storage capacity available in each access hub for the following period. In this study, we assume capacity module relocation is performed by a separate business unit whose routing decisions are out of the scope of the research reported in this paper.  \\
The objective is to minimize the cost incurred by operating such a network of access hubs without disrupting underlying operations. The decision scope is tactical (capacity deployment) and requires the integration of operational decisions. However, since the main interest is a set of tactical decisions, there is no need to explicitly model operations, but only to approximate the impact of deployment decisions on routing cost and time synchronization. \\
Let $L$ be a set of access hub locations and $W$ a set of external depots composing a network $G=(N=L\cup W,A)$ where $A$ is the complete set of directed arcs between locations in $N$. A capacity deployment of $I_0$ capacity modules in time $t$ over the network is represented by a vector $S(t) = (S_l(t), \forall l \in N)$. The relocation of capacity modules can be represented as vectors $R(t)=(R_a(t),\forall a\in A)$. Accordingly, there are ${I_0+|L|-1 \choose I_0}$ possible arrangements of $I_0$ modules over $|L|$ locations. In the case where $I_0 \geq |L|$ and that each location gets at least one module, there are ${I_0-1 \choose |L|-1}$ possible arrangements. In this realistic context, access hub networks are expected to be composed of a high number of locations (i.e., hundreds). Thus, state and action spaces would be significantly large-sized, which results in curse of dimensionality issues (\cite{powell2007approximate}). \\
Moreover, a set of realization scenarios $\omega \in \Omega$ with probability $\phi_{\omega}$ is considered. The number of pickups and deliveries as well as the storage volume requirements are observed hourly and respectively represented as a vectors $\rho^P(\tau, \omega) = (\rho^P_l(\tau,\omega), \forall l \in L)$, $\rho^D(\tau, \omega) = (\rho^D_l(\tau, \omega), \forall l \in L)$ and $D(\tau, \omega)=(D_l(\tau, \omega), \forall l \in L)$, for every operations hour  $\tau \in T_t$, where $t \in T$ is an operations horizon between two deployment periods (e.g. a week). If a courier or rider observes a lack of storage capacity when visiting an access hub, the courier or rider can
perform the following recourse actions: pool capacity by making a detour towards a neighboring access hub
with extra capacity or consign its load to a nearby third-party business (e.g. local shop) for a certain
price agreed upon (uncapacitated recourse). Once volume requirements are observed, recourse actions are taken for each operational period $\tau$: capacity pools as a vector $P(\tau, \omega)=(P_a( \tau,\omega),\forall a\in A_{pool})$ where $A_{pool}$ is the set of arcs on which capacity can be pooled, and consignments as a vector $Z(\tau,\omega)=(Z_l(\tau,\omega),\forall l\in L)$. At any time $\tau$ in scenario $\omega$, the system can thus be represented as a state $S_t = S(t) \text{ s.t. } \tau \in T_{t}$ and an action $x_\tau = (R(\tau), P(\tau,\omega), Z(\tau,\omega)) \text{ s.t. } \tau \in T_{t}$, where $R(\tau)$ is the null vector except for $\tau = t, \forall t \in T$. Based on the optimisation framework  proposed in \cite{powell2019unified}, our stochastic optimization challenge for the access hub dynamic pooled capacity deployment problem can be formulated as follows:
\begin{equation} \label{opt challenge}
    \min_{x_\tau \in X(\tau)} \mathop{\mathbb{E}}_{\omega \in \Omega}\left \{\sum_{t \in T} \sum_{\tau \in T_t} C_\tau(S_t,x_\tau,\rho^P(\tau,\omega), \rho^D(\tau,\omega), D(\tau,\omega))|S_0 \right \}
\end{equation}
where $X(\tau)$ is the set of feasible actions at time $\tau$, $S_0$ is the initial state of the system, and $C_\tau(\cdot)$ is the cost function at time $\tau$. Figure \ref{fig: Timeline} illustrates the dynamics of the problem with the tactical decision timeline: before each period $t$, a network deployment strategy $S(t)$ is decided through relocation decisions $R(t)$ and implemented right before the beginning of period $t$. Then, demand realized and recourse actions are taken in each period $\tau \in T_t$. At the end of periods $T_t$, a network deployment strategy $S(t+1)$ is decided through relocation decisions $R(t+1)$ and implemented right before the beginning of period $t+1$ and the process repeats.

\subsection{Operations Cost Approximation and Synchronization Modeling \label{under op}}

\begin{figure}[H]
\centering
\fbox{\includegraphics[scale=0.30]{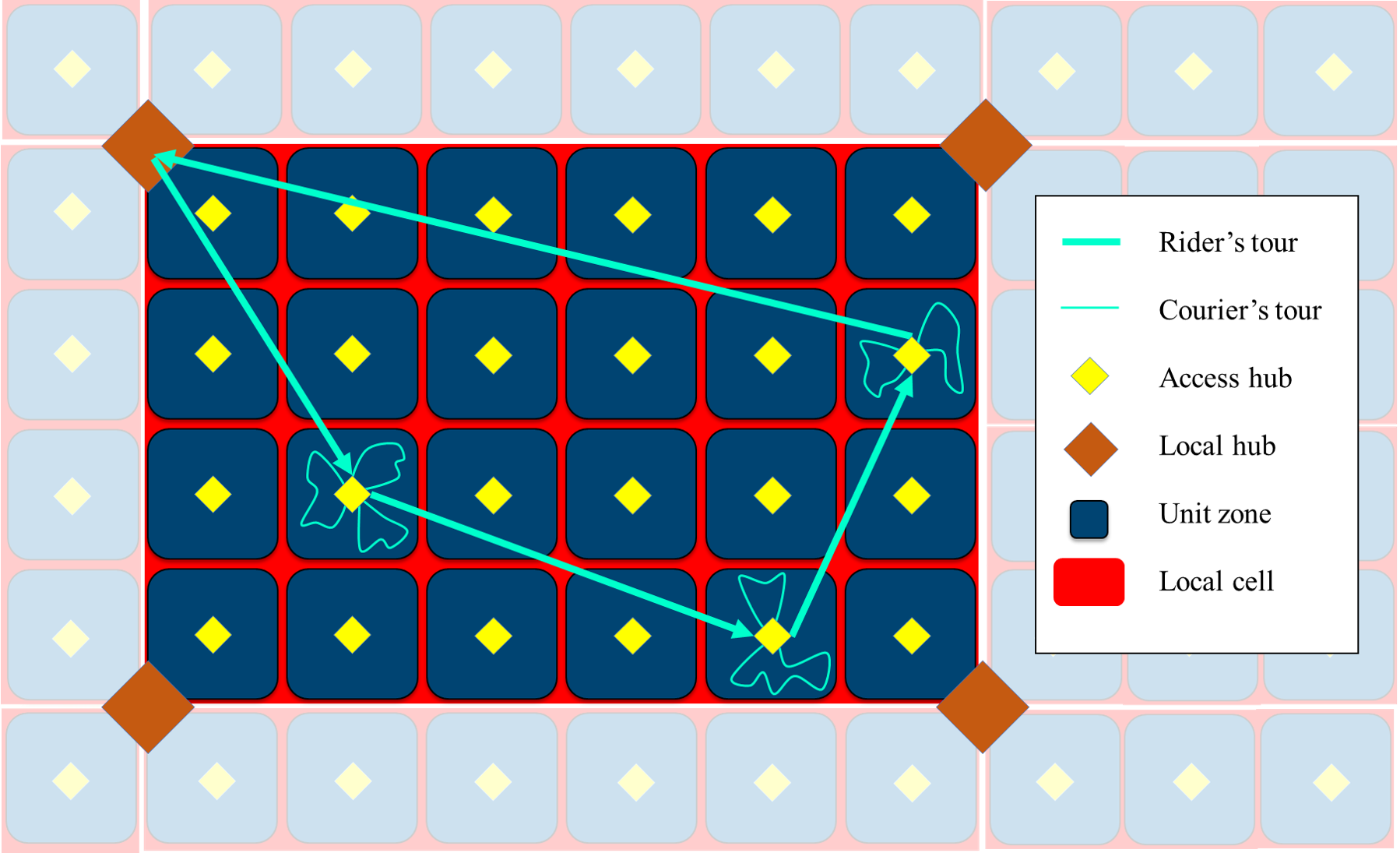}}
\caption{Illustration of Couriers and Riders Operations}
\label{fig: Operations}
\end{figure}

Once decisions on capacity deployment are set for a given period $t$, they strongly impact the quality of operations performed by couriers and riders. More specifically, capacity at each location impacts the number and costs of detour and perturb the synchronisation of the operations between couriers and riders at each location. Accordingly, the surrounding objective of integrating routing operations is to evaluate the performance of the capacity deployment in minimizing the detours due to an  underestimation of the capacity needs and in guaranteeing the synchronisation of the operations between couriers and riders at each location. To do so, this subsection proposes to develop routes with detours cost approximations, and travel time approximations. It builds on a refined granularity of routing operations periods (hourly) and uncertain storage volume requirements.
\\
Figure \ref{fig: Operations} illustrates the use of access hubs in first/last mile parcel logistics operations during a period $\tau$ with one rider and 3 couriers. At the operational level, pickup and delivery decisions are made hourly ($\tau$) based on the volume-based capacity made available at each access hub. In addition, capacity relocation determines the number and costs of recourse actions needed to satisfy the requested volumes. With the consideration of capacity pooling recourse, the pickup and delivery problem with transshipment faced by couriers and riders adds the feature of detours. 
Here, to ensure timely transshipment operations, the detours performed by couriers and riders, are limited to their original time period ($\tau$), avoiding couriers and riders to be desynchronized.
Since these detours necessitate additional moves and are time consuming, this comes with a supplementary incurred cost.

It is clear that capturing the dynamics of underlying operations  when taking capacity deployment decisions leads to better solutions. However, the pickup and delivery problem with transshipment is NP-hard (\cite{RaisPandDwithTransshipment}) and including it explicitly in the tactical model would make it intractable. Since the goal is to foster best capacity deployment decisions, it is sufficient to anticipate the operations costs and time synchronisation constraints using scenario-based continuous approximations.

Accordingly, hereafter is proposed a tractable approximation of each period $\tau$ pickup and delivery problem with transshipment by developing deterministic continuous approximations of vehicle routing problems.
The starting point of the proposed approximations is the estimation
of the vehicle routing problem length when the depot (from which vehicles start their routes) is not necessarily located in the area where customers are located as proposed in \cite{Daganzo94}:
\begin{equation} \label{VRPDaganzo}
VRP(n) = 2rm + nk(\delta)^{-\frac{1}{2}}
\end{equation}
where $r$ is the average distance between the depot and the customer locations, $m$ is the number of routes required to serve all customers, $n$ is the number of points to be visited, $k$ is a constant parameter that can be estimated through simulation (\cite{DaganzoBook}), and $\delta$ is the density of points in the area. An a priori lower bound on the number of routes required to serve all customers, $m$, is $n/Q$ where $Q$ is the capacity of one vehicle in terms of customer locations. The first term of approximation (\ref{VRPDaganzo}) represents the line-haul (back and forth) performed by vehicle to travel from the depot to the area where customers are located, and the second term represents the tour performed by traveling between each successive stops. Based on these seminal works, the next subsection proposes an adaptation of these equations to the operational context of riders and couriers, and develops an explicit time-based estimation of their operations.
\\

\subsubsection*{Riders operations}
\begin{figure}[H]
\centering
\fbox{\includegraphics[scale=0.30]{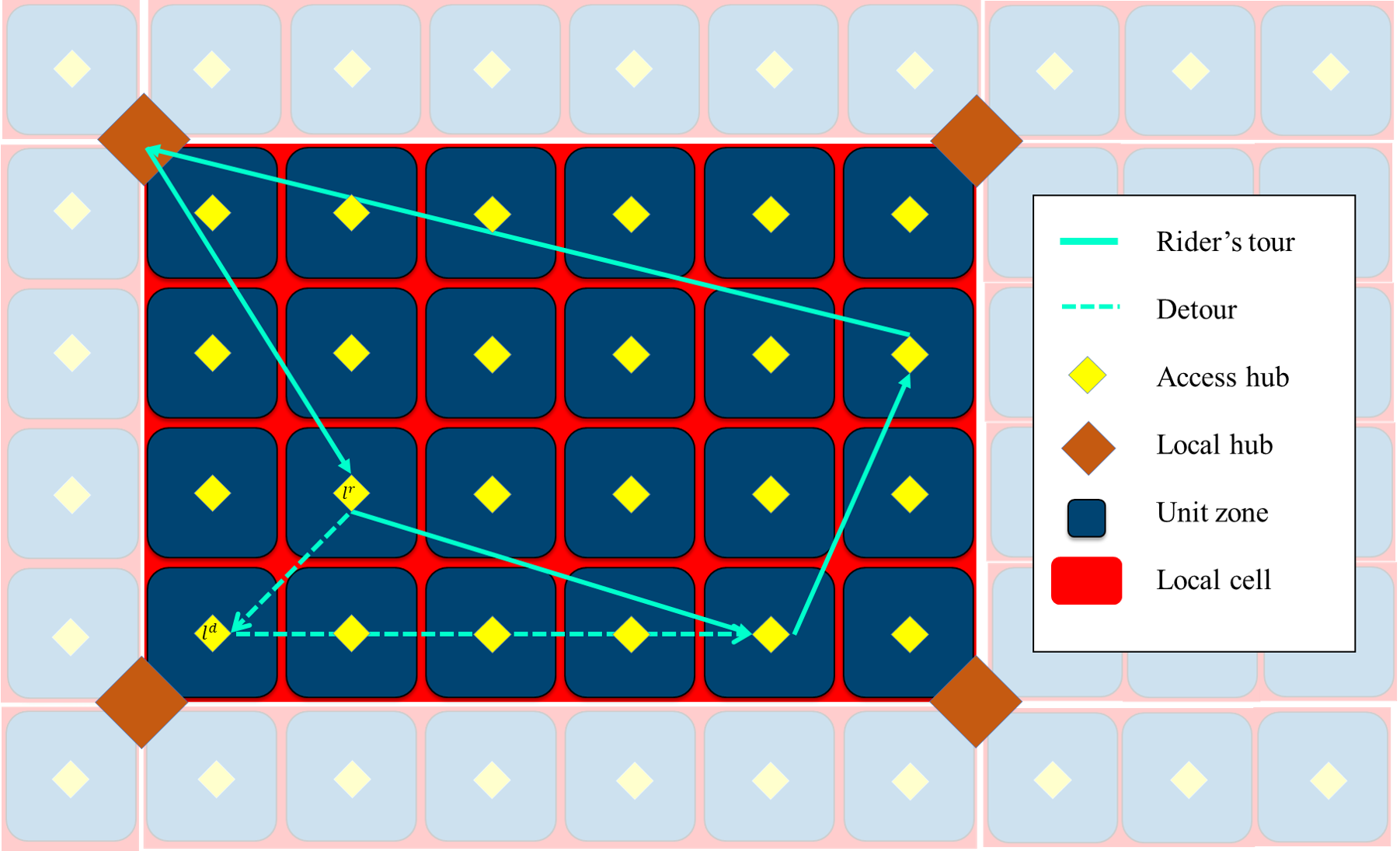}}
\caption{Illustration of a Rider's Tour with Detour}
\label{fig: Rider}
\end{figure}

Riders work in local cells, which are clusters of access hubs served by the same upper level local hub(s) as illustrated in Figure \ref{fig: Rider}. Riders visit a set of $n_{LC}$ access hubs within their local cell of area $A_{LC}$ (and density $\delta_{LC} = \frac{n_{LC}}{A_{LC}}$) to pickup and deliver  parcels as part of a defined route (e.g. planned beforehand based on averaging  network's load). At the time of deployment, underlying  riders' routes are not known with certainty, but need to be estimated in order to anticipate operations performance. When a rider makes his tour in period $\tau$ under scenario $\omega$ two cases are possible: (i) the tour is operated as planned because sufficient capacity is deployed at all visited access hubs in the route or because the detours are assigned to access hubs that are already in the remaining itinerary of the rider (bold lines in Figure \ref{fig: Rider}); (ii) the rider tour is perturbed due to a lack of capacity at an access hub, and thus has to perform an immediate detour to a neighboring access hub before pursuing the rest of the regular tour (dash lines in Figure \ref{fig: Rider}).

Given approximation (\ref{VRPDaganzo}), if the number of detours performed by riders in local cell $LC$ in period $\tau$ in scenario $\omega$ is $n_{LC}^R(\tau, \omega)$, the route length estimation with detours of riders' operations is:

\begin{equation} \label{VRPRider}
VRP^R_{LC}(\tau, \omega) = 2r_{LC}m_{\tau}^R(\omega) + (n_{LC}+n_{LC}^R(\tau, \omega))k^R(\delta_{LC})^{-\frac{1}{2}}
\end{equation}
where $n_{LC}$ is the total number of access hubs in local cell $LC$, $r_{LC}$ is the average distance between $LC$'s local hub(s) and its access hubs, and $m_{\tau}^R(\omega)$ is the number of riders' operating. \\
The cumulative time (in time-rider) necessary to perform tours approximated in (\ref{VRPRider}) is:
\begin{equation} \label{rider time cumulative}
T^R_{LC}(\tau, \omega) = m_{\tau}^R(\omega)(t^R_s + \frac{2r_{LC}}{s^R_0}) + (n_{LC} + n_{LC}^R(\tau, \omega))(\frac{ k^R(\delta_{LC})^{-\frac{1}{2}}}{s^R} + t^R_a) + (\sum_{l \in LC}(\rho_l^D(\tau,\omega) + \rho_{l}^P(\tau,\omega)))t^R_u
\end{equation}
where the first term is the time spent to setup tours ($t_s^R$ per tour) and perform the line-haul at a speed of $s_0^R$, the second term represent the travel time between stops at a speed of $s^R$ and the stopping time $t_a^R$ per access hub, and the third term represents the service time (handling) $t_u^R$ per pickup and delivery. \\
Thus, the cost associated with riders' operations in local cell $LC$ in period $\tau$ in scenario $\omega$ is:
\begin{equation}
C^R_{LC}(\tau, \omega) = m_{\tau}^R(\omega)(c^R_f + 2r_{LC} c^R_{v_0}) + (n_{LC} + n_{LC}^R(\tau, \omega))k^R(\delta_{LC})^{-\frac{1}{2}} c^R_v + T^R_{LC}(\tau, \omega) c^R_w
\end{equation}
where the first term represents the fixed, $c^R_f$, and variable, $c^R_{v_0}$ in line-haul and $c^R_v$ in tour, costs associated with vehicles, and the second term represents the variable labor cost $c^R_w$ of $m_\tau^R(\omega)$ riders. \\
Since the nominal routing cost (with no detours) is a sunk cost incurred regardless of the capacity deployment, the marginal cost is sufficient to inform the tactical decision of the impact of recourse actions. The marginal cost of the detours induced by the tactical decisions, or difference between the rider routing cost with detours and the nominal rider routing cost, is:
\begin{equation}
\Delta C^R_{LC}(\tau, \omega) = n_{LC}^R(\tau, \omega)k^R(\delta_{LC})^{-\frac{1}{2}} c^R_v + \Delta T^R_{LC}(\tau, \omega) c^R_w
\end{equation}
where the time associated with performing detours is the time needed to perform detours:
\begin{equation}
    \Delta T^R_{LC}(\tau, \omega) = n_{LC}^R(\tau, \omega)(\frac{ k^R(\delta_{LC})^{-\frac{1}{2}}}{s^R} + t^R_a)
\end{equation}

\subsubsection*{Couriers operations}

\begin{figure}[H]
\centering
\fbox{\includegraphics[scale=0.35]{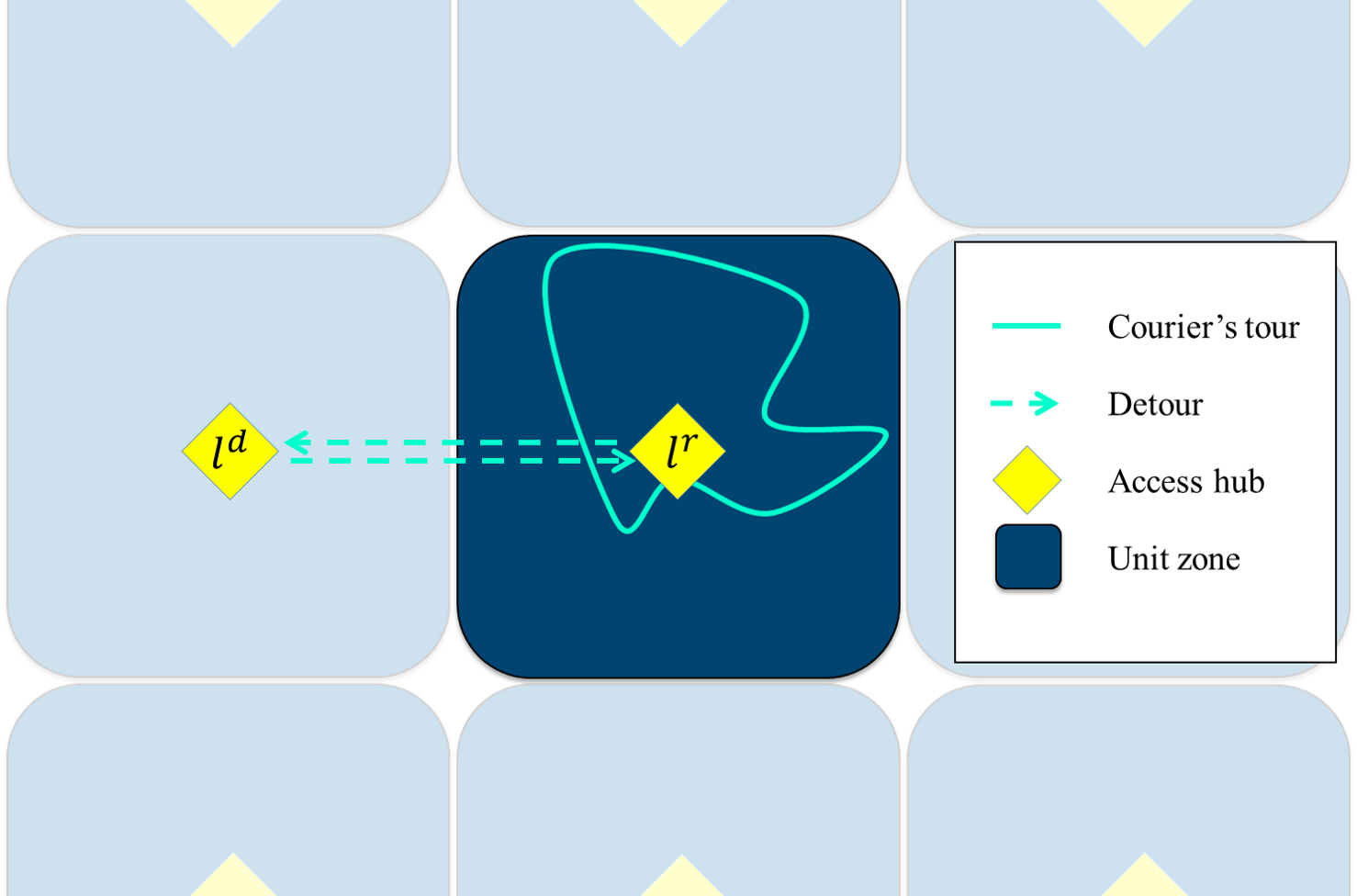}}
\caption{Illustration of a Courier's Tour with Detour}
\label{fig: Courier}
\end{figure}

Couriers operate in unit zones, which are clusters of pickup and delivery points served by access hub(s). Couriers leave their reference access hubs to visit customers and perform pickups/deliveries before returning to their access hub. When a courier arrives at the courier’s access hub with picked parcels, if the courier observes a lack of capacity, the courier can be immediately directed to available capacity in some neighboring access hub. Then, the courier will perform a detour (out and back) to the assigned neighbour access hub before starting their next tour from their reference access hub. Figure \ref{fig: Courier} illustrates a courier's tour and a detour as described.

Since access hubs are located in the same area as pickup/delivery locations, the line-haul distance at this echelon is negligible, which eliminates the first term of approximation (\ref{VRPDaganzo}). If the number of detours performed by couriers on arc $a \in A_{pool}(l)= \{a=(l,j), \forall j: (l,j) \in A_{pool}\}$ of length $d_a$ in period $\tau$ under scenario $\omega$ is $n_{a}^C(\tau, \omega)$, the route length estimation with detours of couriers' operations is:
\begin{equation} \label{VRPCouriers}
VRP^C_{l}(\tau, \omega) = (\rho_{l}^P(\tau,\omega) + \rho_{l}^D(\tau,\omega))k^C(\delta_{l})^{-\frac{1}{2}} + \sum_{a \in A_{pool}(l)} (2 n_{a}^C(\tau, \omega) d_a)
\end{equation}
where the first term represents the total length of tours performed by couriers to visit pickup/delivery locations, and the second term represents the detours (out and back) performed between access hub $l$ and its neighboring access hubs. \\
The cumulative time (in time-courier) necessary to perform courier tours is based on the approximation in (\ref{VRPCouriers}) as follows:
\begin{equation} \label{courier time cumulative}
T^C_{l}(\tau, \omega) = (\rho_{l}^P(\tau,\omega) + \rho_{l}^D(\tau,\omega))(\frac{ k^C(\delta_{l})^{-\frac{1}{2}}}{s^C} + t^C_a) + \sum_{a \in A_{pool}((l)} n_{a}^C(\tau, \omega)(\frac{2 d_a}{s^C_0} + t^C_a) + (\rho_{l}^D(\tau,\omega) + \rho_{l}^P(\tau,\omega))t^C_u
\end{equation}
where the first term represents the travel time between pickup/delivery locations at a speed of $s^C$ and the stopping time $t_a^C$ per stop, the second term represents the travel time during detours to neighboring access hubs at a speed of $s^C_0$ plus a stopping time $t_a^C$, and the third term represents the service time (handling) $t_u^C$ per pickup and delivery. \\
Thus, the cost associated with couriers' operations at access hub $l$ in period $\tau$ under scenario $\omega$ is:
\begin{equation}
C^C_{l}(\tau, \omega) = ((\rho_{l}^P(\tau,\omega) + \rho_{l}^D(\tau,\omega))k^C(\delta_{l})^{-\frac{1}{2}} c^C_v + \sum_{a \in A_{pool}(l)} (2 n_{a}^C(\tau, \omega) d_a) c^C_{v_0}) + T^C_{l}(\tau, \omega) c^C_w
\end{equation}
where the first term represents the variable travel costs, respectively $c^C_v$ between pickup/delivery locations and $c^C_{v_0}$ between access hubs, and the second term represents the variable labor cost $c^C_w$ of $m_\tau^C(\omega)$ couriers. \\
Again, since the nominal routing cost (with no detours) is a sunk cost incurred regardless of the capacity deployment, the marginal cost is sufficient to inform the tactical decision of the impact of recourse actions.The marginal cost of the detours induced by the tactical decisions, or the difference between the courier routing cost with detours and the nominal courier routing cost is:
\begin{equation}
\Delta C^C_{l}(\tau, \omega) = \sum_{a \in A_{pool}(l)} (2 n_{a}^C(\tau, \omega) d_a) c^C_{v_0}) + \Delta T^C_{l}(\tau, \omega) c^C_w
\end{equation}
where the time associated with performing detours is:
\begin{equation}
    \Delta T^C_{l}(\tau, \omega) = \sum_{a \in A_{pool}(l)} n_{a}^C(\tau, \omega)(\frac{2 d_a}{s^C_0} + t^C_a)
\end{equation}

\subsubsection*{Operations Synchronization}
Recall that a key objective of integrating routing operations with the capacity deployment problem is to  guarantee the synchronisation of the operations between couriers and riders at each location. To do so, this subsection proposes to develop time-based synchronisation constraints based on the travel time approximations (\ref{rider time cumulative}) and (\ref{courier time cumulative}), developed above.

Parcels transshipped from riders to couriers and couriers to riders through access hubs must be transshipped during the period of time the parcels are within the network. That is, the length of a courier's (respectively rider's) original tour, plus the added detour(s) must not exceed the maximum length feasible within one operational period.
For riders’ operations, at the local cell level, this tour length can be expressed, based on the number of riders ($m_\tau^R(\omega$) in period $\tau$ under scenario $\omega$, as follows:
\begin{equation} \label{Rider_Synchro_Constraint}
    T^R_{LC}(\tau, \omega) \leq m_\tau^R(\omega) \Delta_\tau, \forall \omega \in \Omega, LC \in \mathcal{LC}, \tau \in T_t, t \in T
\end{equation}
where $\Delta_\tau$ is the length of period $\tau$.
Similarly, for couriers' operation, at the access hub level, synchronization can be expressed, based on the number of couriers ($m_\tau^C(\omega$) in period $\tau$ under scenario $\omega$, as follows:
\begin{equation} \label{Courier_Synchro_Constraint}
    T^C_{l}(\tau, \omega) \leq m_\tau^C(\omega) \Delta_\tau, \forall \omega \in \Omega, l \in L, \tau \in T_t, t \in T
\end{equation}

\subsection{Two-Stage Stochastic Program Formulation for the Access Hub Dynamic Pooled Capacity Deployment Problem  }
In this section, a stochastic programming formulation is proposed to tackle the optimization problem (\ref{opt challenge}) presented in section \ref{biz}.
We remark that the stochastic optimisation problem (\ref{opt challenge}) can be modeled as a multi-stage stochastic program based on a scenarios tree. However, this program would be intractable for realistic size instances, due to its combinatorial structure and non-anticipatory constraints \citep{Schultz2003}.
Under a rolling horizon framework, the model is built here on the relaxation approach \citep{Shapiro2009} that is applied to transform the multi-stage stochastic program to a two-stage stochastic program with multiple tactical periods.
More specifically, it consists in transferring all the capacity deployment decisions of the $T$ periods to the first-stage in order to be set at the beginning of the horizon. In this case, only first-stage design decisions ($t=1$) are made here and now, but subsequent capacity deployment decisions ($t>1$) are deferrable in time according to their deployment period.
Hereafter are introduced the additional sets, input parameters, random variables and decision variables that formulate the overall model.

\subsubsection*{Sets}
\begin{align*}
\begin{tabularx}{\textwidth}{ll}
$\mathcal{L}$ & access hub locations, indexed by $l$ \\
$\mathcal{LC}$ & local cells, indexed by $LC$ \\
$\mathcal{W}$ & depot locations, indexed by $l$ \\
$\mathcal{A}$ & arcs between two locations of the network $\mathcal{L}\cup\mathcal{W}$, indexed by $a$ \\
$\mathcal{G}$ & asymmetric graph $(\mathcal{L}\cup\mathcal{W},A)$ satisfying the triangle inequality \\
$T$ & tactical periods, indexed by $t$, covering the planning horizon\\
$T_t$ & subset of operational (demand) periods, indexed by $\tau$, between periods $t$ and $t+1$ \\
$\Omega$ & scenarios, indexed by $\omega$ \\
$\delta^+(l)$ & incoming relocation arcs in location $l \in \mathcal{G}$ \\
$\delta^-(l)$ & outgoing relocation arcs from location $l \in \mathcal{G}$ \\
$N^+(l)$ & incoming recourse arcs in location $l$ ; $N^+(l) \subset \delta^+(l)$\\
$N^-(l)$ & outgoing recourse arcs from location $l$ ; $N^-(l) \subset \delta^-(l)$\\
$A_{pool}(l)$ & recourse arcs available for capacity pooling from location $l$
 \end{tabularx}
 \end{align*}
\\
\subsubsection*{Input Parameters}

\begin{align*}
\begin{tabularx}{\textwidth}{ll}
$h_l$ & cost of holding one capacity module at location $l$. \\ 
$I_0$ & total number of capacity modules available in the system\\
$\phi_\omega$ & probability of scenario $\omega$\\
$p_l$ & penalty for lacking capacity in location $l$ \\
$r_a$ & cost of relocating one capacity module on $a$. \\ 
$\overline{S_l}$ & maximum number of capacity modules that can be placed in location $l$ \\
$v$ & volume provided by a capacity module \\
$v^R$ & volume that a rider can carry on a tour \\
$v^C$ & volume that a courier can carry on a tour
\end{tabularx}
 \end{align*}
\\
\subsubsection*{Random Variables}

\begin{align*}
\begin{tabularx}{\textwidth}{ll}
$D_l(\tau, \omega)$ & volume requirements in location $l$ in scenario $\omega$ in period $(\tau)$
\end{tabularx}
 \end{align*}
\\
\subsubsection*{Decision Variables}

\begin{align*}
\begin{tabularx}{\textwidth}{ll}
$S_l(t)$ & number of capacity modules available in location $l$ for period $t$ \\
$R_a(t)$ & number of capacity modules relocated through arc $a$ at the beginning of period $t$ \\
$P_a(\tau, \omega)$ & volume shared from location $i$ to $j$, $a=(i,j) \in N_{\omega, \tau}^-(i)$ in period $\tau$ under scenario $\omega$ \\
$Z_l(\tau, \omega)$ &  lack of capacity in volume at location $l$ in period $\tau$ under scenario $\omega$ \\
$n^R_a(\tau, \omega)$ &  number of detours performed by riders on arc $a$ in period $\tau$ under scenario $\omega$ \\
$n^R_{LC}(\tau, \omega)$ &  number of detours performed by riders in local cell $LC$ in period $\tau$ under scenario $\omega$ \\
$n^C_a(\tau, \omega)$ &  number of detours performed by couriers on arc $a$ in period $\tau$ under scenario $\omega$ \\
$n^C_l(\tau, \omega)$ &  number of detours performed by couriers from location $l$ in period $\tau$ under scenario $\omega$
\end{tabularx}
 \end{align*}

\subsubsection*{Model}

\begin{align}
\min &  \sum_{t \in T} \bigg( \sum_{l \in \mathcal{L}\cup \mathcal{W}} h_l S_l(t) + \sum_{a \in A} r_a R_a(t) \nonumber \\
& + \sum_{\omega \in \Omega} \phi_{\omega} \bigg( \sum_{\tau \in T_{t}} \bigg( \sum_{l \in L} (\Delta C_l^C(\tau, \omega) + p_l Z_l(\tau, \omega)) + \sum_{LC \in \mathcal{LC}} \Delta C_{LC}^R(\tau, \omega)  \bigg) \bigg) \bigg) \label{Obj} \\
\text{s.t.: \ } \nonumber \\
& \text{Inventory balance of capacity modules at all locations:} \nonumber \\
& S_l(t) = S_l(t-1) + \sum_{a \in \delta^+(l)} R_a(t) - \sum_{a \in \delta^-(l)} R_a(t), \forall l \in \mathcal{L}\cup \mathcal{W}, t \in T \label{IB} \\
& \text{Total capacity module inventory constraint:} \nonumber \\
& \sum_{l \in \mathcal{L}\cup\mathcal{W}}S_l(t) = I_0, \forall t \in T \label{Inv} \\
& \text{Spatial constraint at all locations:} \nonumber \\
& S_l(t) \leq \overline{S_l}, \forall l \in \mathcal{L}, t \in T \label{Spatial} \\
& \text{Volume requirements satisfaction constraints:} \nonumber \\
& v S_l(t) + \sum_{a \in N^+(l)} P_a (\tau, \omega) - \sum_{a \in N^-(l)} P_a (\tau, \omega) + Z_l (\tau, \omega) \geq D_l(\tau, \omega), \forall l \in L, \tau \in T_{t}, t \in T, \omega \in \Omega \label{Vol Req}\\
& \text{Synchronization constraint for riders' operations: (\ref{Rider_Synchro_Constraint})} \nonumber \\
& \text{Synchronization constraint for couriers' operations: (\ref{Courier_Synchro_Constraint})} \nonumber \\
& \text{Rider's detours count:} \nonumber \\
& n_{LC}^R(\tau, \omega) \geq \sum_{l \in LC} \sum_{a \in A_{pool(l)}} n_a^R(\tau, \omega), \forall LC \in \mathcal{LC}, \tau \in T_{t}, t \in T, \omega \in \Omega \\\label{LC detour count}
& n_a^R(\tau, \omega) \geq \frac{P_a(\tau, \omega)}{v^R}, \forall a \in A_{pool}(l), l \in L, \tau \in T_{t}, t \in T, \omega \in \Omega \\\label{Rider detour count}
& \text{Courier's detours count:} \nonumber \\
& n_{l}^C(\tau, \omega) \geq \sum_{a \in A_{pool}(l)} n_a^C(\tau, \omega), \forall l \in L, \tau \in T_{t}, t \in T, \omega \in \Omega \\\label{AH detour count}
& n_a^C(\tau, \omega) \geq \frac{P_a(\tau, \omega)}{v^C}, \forall a \in A_{pool}(l), l \in L, \tau \in T_{t}, t \in T, \omega \in \Omega \\\label{Courier detour count}
& \text{Integrality and non-negativity constraints:} \nonumber \\
& P_a(\tau, \omega), Z_l(\tau, \omega), n_a^C(\tau, \omega), n_a^R(\tau, \omega), n_{l}^C(\tau, \omega)  \geq 0  \\
& S_l(t), R_a(t)\text{ integer} \label{integrality constraint}
\end{align}


Minimizing expression (\ref{Obj}) corresponds to minimizing the last-mile cost, defined in this paper as the cost of deploying capacity modules in each access hub locations (holding costs) and the relocation costs for each capacity module movement for each reconfiguration period, and the marginal cost incurred by recourse actions (capacity pool from neighboring location and consignment). Constraints (\ref{IB}) and (\ref{Inv}) enforce the conservation of the total number of capacity modules in the network. Constraints (\ref{Spatial}) limit the number of capacity modules that can be deployed in each access hub locations. Constraints (\ref{Vol Req}) enforce that all demand in terms of volume requirement is served by a combination of capacity modules, capacity pools and consignments, in each demand period of each scenario. Constraints (\ref{Rider_Synchro_Constraint}) and (\ref{Courier_Synchro_Constraint}) are the synchronization constraints for the underlying riders and couriers problems as developed in section (\ref{under op}). Constraints (\ref{LC detour count}) and (\ref{Rider detour count}) count the number of detours performed by riders within each local cell based on recourse capacity pooling decisions and the carrying capacity of riders. Constraints (\ref{AH detour count}) and (\ref{Courier detour count}) count the number of detours performed by couriers from each access hub based on recourse capacity pooling decisions and the carrying capacity of couriers.

\section{Solution Approach} \label{SolutionApproach}

In this section, our rolling horizon solution approach is presented, which builds on solving sequentially the two-stage model presented above using scenario sampling, Benders decomposition and acceleration methods. It approximates optimization problem (\ref{opt challenge}) by planning for one capacity deployment period, $t$, at the time and deferring subsequent capacity deployment decisions to the following iterations of the Algorithm. In order to enhance the quality of the solutions produced at each iteration, a $\theta$ tactical lookahead is considered to plan for $1 + \theta$ tactical periods, where only the first period is implementable and the subsequent ones are used as an evaluation mechanism.
The proposed rolling horizon solution approach is described in Algorithm 1.
\noindent
Here, the length of the sub-horizon is controllable; it can represent one tactical period (i.e. myopic, $\theta = 0$) or several of them (i.e. lookahead, $\theta \geq  1$). Of course, when dealing with large-scale networks, the selection of the lookahead length is part of the trade-offs necessary to make in order to keep the model tractable. In order to enhance the solvability of the optimization model (\ref{Obj}-\ref{integrality constraint}), for each sub-horizon $[t, t+\theta]$, a tailored Benders decomposition approach is developed, that fits with the two-stage and multi-period setting of our formulation. It is applied under a large sample of multi-period scenarios. The following subsections address the decomposition approach as well as the associated acceleration methods developed.\\

\begin{algorithm}[h]
\SetAlgoLined
\KwResult{$S_l(t), R_a(t)$}
$S_l(t_0) \longleftarrow S_l(t_0)$\;
\For{$t \in T$}{
$S_l(t), R_a(t) \longleftarrow $ Optimal solution of (\ref{Obj}-\ref{integrality constraint}) for sub-horizon $[t, t+\theta]$;
}
\caption{Rolling Horizon Algorithm with Tactical Lookahead}
\label{dynamic algorithm}
\end{algorithm}

\subsection{Benders Decomposition} \label{BD}

Benders decomposition is a row generation solution method for solving large scale optimization problems by partitioning the decision variables in first stage and second stage variables (\cite{Benders2005}). The model is first projected onto the subspace defined by the first stage variables, replacing the second stage variables by an incumbent; the resulting model is called the restriced master problem. Then, a linear problem with the second stage variables and a candidate solution from the restricted master problem is formulated; the resulting model is called the subproblem and can often be decomposed in independent subproblems. From the solution of the subproblem, feasibility and optimality cuts can be identified and added to the restricted master problem. The algorithm terminates when the incumbent in the restriced master problem is equal to the the value of the subproblem. \\
Suppose the capacity deployment and relocation decisions (first stage decision variables) $S_l(t)$, $S_l(t+1)$, ... ,$S_l(t+\theta)$ and $R_a(t)$, $R_a(t+1)$, ..., $R_a(t+\theta)$ are given with values $\widehat{S_l}(t)$, $\widehat{S_l}(t+1)$, ... , $\widehat{S_l}(t+\theta)$ and $\widehat{R_a}(t)$, $\widehat{R_a}(t+1)$, ..., $\widehat{R_a}(t+\theta)$. Then, the subproblem can be defined as taking recourse action decisions (i.e. second stage decisions; capacity pooling) to minimize the approximate overall operations costs. The subproblem can be decomposed per scenario $\omega$, operational period $\tau$ and local cell $LC$ into a set of independent subproblems as follows:
\begin{align} \label{subproblem}
 SP_{LC}(\tau, \omega) =  \min &  \sum_{l \in L(LC)} (\Delta C_l^C(\tau, \omega) + p_l Z_l(\tau, \omega)) + \Delta C_{LC}^R(\tau, \omega) \\
\text{s.t.: \ } \nonumber \\
& \text{Volume requirements satisfaction constraints:} \nonumber \\
& v \widehat{S_l}(t) + \sum_{a \in N^+(l)} P_a(\tau, \omega) - \sum_{a \in N^-(l)} P_a(\tau, \omega) + Z_l (\tau, \omega) \geq D_l(\tau, \omega), \forall l \in L(LC) \label{dual1}\\
& \text{Synchronization constraint for riders' operations: (\ref{Rider_Synchro_Constraint})}  \nonumber \\
& \text{Synchronization constraint for couriers' operations: (\ref{Courier_Synchro_Constraint})} \nonumber \\
& \text{Detour linking constraints: (\ref{LC detour count}), (\ref{Rider detour count}), (\ref{AH detour count}), (\ref{Courier detour count})} \nonumber \\
& P_a(\tau, \omega), Z_l(\tau, \omega), n_a^C(\tau, \omega), n_a^R(\tau, \omega), n_{l}^C(\tau, \omega) \geq 0 \nonumber
\end{align}


It is important to notice that the defined subproblems are feasible regardless of the value of the tactical decisions (first stage variables); This is possible thanks to the variables $Z_l(\tau, \omega)$ that compensate for any lack of capacity in the network by incurring a large cost. \\
Solving each subproblem using a dualization strategy, one can identify the following optimality cuts for each local cell, operational period $\tau$ and scenario $\omega$:

\begin{flalign} \label{Benders cuts}
    q_{LC}(\tau, \omega) \geq & \sum_{l \in L(LC)}\pi^j_l(\tau,\omega)(D_l(\tau, \omega)-v S_l(t)) \\\nonumber
    &+ \mu^j_{LC}(\tau, \omega) \bigg( m^R(\tau, \omega)(\Delta_\tau - (t^R_s + \frac{2r_{LC}}{s^R_0})) - n_{LC}(\frac{ k^R(\delta_{LC})^{-\frac{1}{2}}}{s^R} + t^R_a)  \\ \nonumber
    & \qquad \qquad \qquad -\sum_{l \in LC}(\rho_l^D(\tau, \omega) +  \rho_l^P(\tau, \omega))t^R_u  \bigg) && \\ \nonumber
    &+ \sum_{l \in L(LC)} \bigg( \lambda^j_l(\tau,\omega)( m^C(\tau,\omega) \Delta_\tau - (\rho_l^P(\tau, \omega) + \rho_l^P(\tau, \omega))(\frac{ k^C(\delta_{l})^{-\frac{1}{2}}}{s^C} + t^C_a + t^C_u) ) \bigg) &&
\end{flalign}

\noindent where $j \in J$, the set of extreme points of the dualized subproblem; $\pi^j_l(\tau,\omega)$, $\mu^j_{LC}(\tau, \omega)$ and $\lambda^j_l(\tau,\omega)$ are the dual values respectively associated with constraints (\ref{dual1}), (\ref{Rider_Synchro_Constraint}) and (\ref{Courier_Synchro_Constraint}).


Finally, the restricted master problem, whose objective minimizes the cost of deploying capacity modules in each access hub and the relocation costs for each capacity module for each period subject to the optimality cuts, can be formulated as follows:

\begin{align}
RMP = \min & \sum_{t}^{t+\theta} \bigg( \sum_{l \in \mathcal{L}\cup \mathcal{W}} h_l S_l(t) + \sum_{a \in A} r_a R_a(t) + \sum_{\omega \in =\Omega} \phi_{\omega} \sum_{\tau \in T_{t}} \sum_{LC \in \mathcal{LC}} q_{LC}(\tau, \omega) \bigg)\\
\text{s.t.: \ } \nonumber \\
& \text{Inventory balance of capacity modules at all locations: (\ref{IB})} \nonumber \\
& \text{Total capacity module inventory constraint: (\ref{Inv})} \nonumber \\
& \text{Spatial constraint at all locations: (\ref{Spatial})} \nonumber \\
& \text{Optimality cuts: (\ref{Benders cuts})}, \forall j \in \overline{J}\subset J\\\nonumber
& S_l(t), R_a(t) \text{ integer} \nonumber
\end{align}
 \noindent Solving the restriced master problem with added optimality cuts provides new values $\widehat{S_l}(t)$ and $\widehat{R_a}(t)$, and a new incumbent solution. This process can be executed iteratively until the incumbent solution equals the subproblem value, indicating optimality.


\subsection{Acceleration Methods} \label{preprocessing and acceleration}

The following subsection describes acceleration methods developed to improve the performance of the proposed solution approach on large instances. The acceleration techniques retained are those that improve significantly the convergence speed of the benders decomposition algorithm for the proposed model.

\subsubsection*{Pareto-optimal Cuts}

The proposed implementation of the benders decomposition can be improved using Pareto-optimal cuts, which requires to solve two linear programs: the original subproblem (\ref{subproblem}), and the Pareto subproblem. The result is the identification of the strongest cut when the original subproblem solution has multiple solutions. A Pareto-optimal solution produces the maximum value at a core point, which is required to be in the relative interior of the convex hull of the subregion defined by the first stage variables. The Pareto subproblem can be decomposed per scenario $\omega$, operational period $\tau$ and local cell $LC$ in a set of independent Pareto subproblems as follows:
\begin{align} \label{Pareto subproblem}
 \min &  \sum_{l \in L(LC)} (\Delta C_l^C(\tau, \omega) + p_l Z_l(\tau, \omega)) + \Delta C_{LC}^R(\tau, \omega) + v_{SP} Y \\
\text{s.t.: \ } \nonumber \\
& v (S_l^0(t) + \sum_{a \in N^+(l)} P_a(\tau, \omega) - \sum_{a \in N^-(l)} P_a(\tau, \omega) + Z_l (\tau, \omega) + (D_l(\tau, \omega) - v \widehat{S_l}(t)) Y \nonumber \\
& \qquad \qquad \qquad \qquad \qquad \qquad \qquad \qquad \qquad \qquad \qquad \qquad \qquad \qquad \qquad \geq D_l(\tau, \omega), \forall l \in L(LC) \label{pareto1}\\
& \text{Modified synchronization constraint for riders' operations:}  \nonumber \\
& T^R_{LC}(\tau, \omega) \leq m^R(\tau, \omega) \Delta_\tau (1 - Y) \nonumber \\
& \qquad \qquad + \bigg( m_{\tau}^R(\omega) \bigg(t^R_s + \frac{2r_{LC}}{s^R_0} \bigg) + n_{LC} \bigg( \frac{ k^R(\delta_{LC})^{-\frac{1}{2}}}{s^R} + t^R_a \bigg) + \bigg( \sum_{l \in LC}(\rho_{l \tau}^D(\omega) + \rho_{l \tau}^P(\omega)) \bigg)t^R_u \bigg) Y \label{pareto_rider}\\
& \text{Modified synchronization constraint for couriers' operations:} \nonumber \\
& T^C_{l}(\tau, \omega) \leq m^C(\tau, \omega) \Delta_\tau (1 - Y) + \bigg( (\rho_{l \tau}^P(\omega) + \rho_{l \tau}^P(\omega)) \bigg(\frac{ k^C(\delta_{l})^{-\frac{1}{2}}}{s^C} + t^C_a + t^C_u \bigg) \bigg) Y, \forall l \in L(LC) \label{pareto_courier}\\
& \text{Detour linking constraints: (\ref{LC detour count}), (\ref{Rider detour count}), (\ref{AH detour count}), (\ref{Courier detour count})} \nonumber \\
& P_a(\tau, \omega), Z_l(\tau, \omega), n_a^C(\tau, \omega), n_a^R(\tau, \omega), n_{l}^C(\tau, \omega), Y \geq 0 \nonumber
\end{align}
\noindent where $v_{SP}$ is the value of the corresponding original subproblem and $S_l^0(t)$ a core point of the current solution to the restricted master problem.  Solving each Pareto subproblem using a dualization strategy, one can identify strengthened optimality cuts (\ref{Benders cuts}) by assigning $\pi^j_l(\tau,\omega)$, $\mu^j_{LC}(\tau, \omega)$ and  $\lambda^j_l(\tau,\omega)$ the dual values respectively associated with constraints (\ref{pareto1}), (\ref{pareto_rider}) and (\ref{pareto_courier}). \\
The proposed implementation also updates the core point, which can be seen as an intensification procedure: locations that are rarely given capacity modules decay toward low values while locations with consistent capacity module presence in every solution are assigned a high coefficient in Pareto solutions. The update rule was introduced in \cite{papadakos2009integrated}, and consists of updating the core point at iteration $k$, $S^{0(k)}$ by combining it with the solution of the master problem at this iteration, $\widehat{S}^{(k)}$, using a factor $\lambda$. \cite{maheo2019benders} suggest that a factor $\lambda=1/2$ yields the best results.  The update rule is defined as follows:

\begin{align*}
    S^{0(k+1)}_l(t) = \frac{S^{0(k)}_l(t) + \widehat{S_l}^{(k)}(t)}{2}, \forall l \in L, t \in T \\
    S^{0(k+1)}_l(t) = \frac{I_0 - \sum_{l' \in L}S^{0(k+1)}_{l'}(t)}{|W|}, \forall l \in W
\end{align*}

\noindent where $k$ is the current iteration of the Benders algorithm.

\subsubsection*{$\epsilon$-optimal Method}
 When dealing with large-scale instances, the $\epsilon$-optimal method as described in \cite{rahmaniani2017benders} has proven to speed up the proposed Benders decomposition algorithms by avoiding to solve the restricted master problem to optimality at each iteration, while guaranteeing an optimal gap within $\epsilon$. It is not necessary to solve the restricted master problem to optimality at each iteration to generate good quality cuts, and there is no incentive to do so at the beginning of the algorithm because the relaxation is weak. Instead, the restricted master problem can be solved with a relaxed optimality gap by adding a constraint forcing the objective value to be improved by at least $\widehat{\epsilon}$ percent compared to the previous solution. Then, when no feasible solution is found, $\widehat{\epsilon}$ is decreased. The same mechanism is applied until $\epsilon$ is reached; the algorithm terminates when no feasible solution is found to the restricted master problem, guaranteeing that the current solution is within $\epsilon$ of the optimal.

\section{Experimental Results} \label{Results}

In this section, the results of numerical experiments are presented in order to validate the developed modeling and solution approaches, and to analyze the performance of the proposed capacity deployment strategy for urban parcel logistics. After describing the test instances which are inspired from the real data of a large parcel express carrier, experimental results about the computational performance of the solution approach are presented. Then, the performance of the dynamic pooled capacity deployment strategy is exposed and compared to its static counterpart. Finally sensitivity analyses are conducted on the capacity pooling distance and the holding costs to derive further insights.

\subsection{Experimental setting} \label{Exp setting description}

\begin{figure}[H]
\fbox{\includegraphics[width=\linewidth]{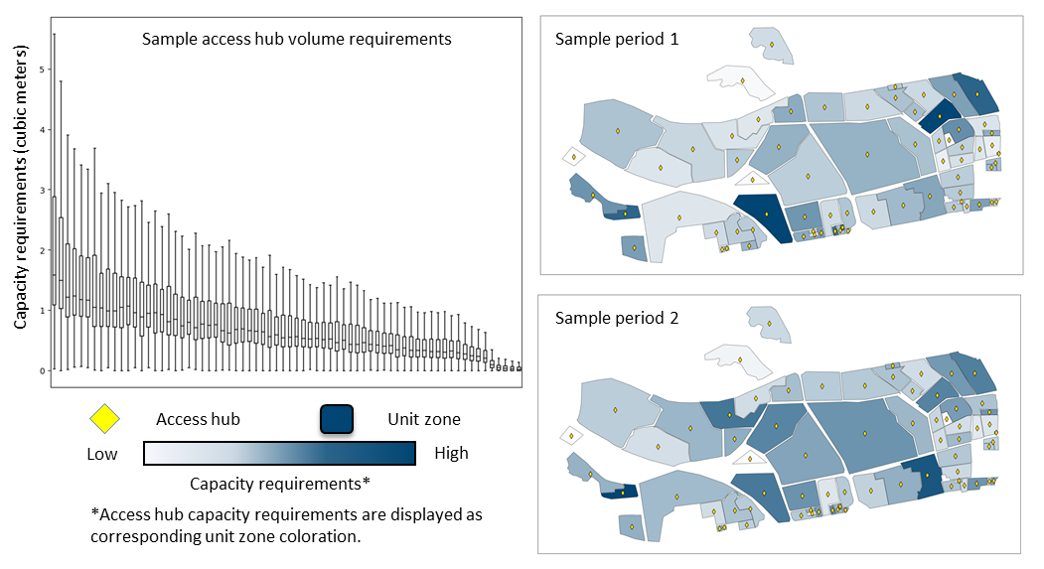}}
\caption{Demand Dynamics Sample for a Local Cell}
\label{demand sample}
\centering
\end{figure}

Table \ref{Instances} summarizes the characteristics of the considered instances: number of access hub locations, number of local cells, and area and population covered by the network. 100 non-stationary demand scenarios are generated randomly from given distributions at the hourly level with monthly, weekly, daily and hourly seasonality factors. Figure \ref{demand sample} illustrates demand dynamics by displaying access hub volume requirements box plots and snapshots of demand levels in two consecutive tactical periods as seen in Figure \ref{fig: Relocation} for a sample local cell from instance E. The number of scenarios is chosen to ensure tactical decision stability with a reasonable in-sample statistical gap  ($1.5\%$) and coefficient of variation ($0.5\%$) as detailed in \ref{Experiment Parameters}. The considered planning horizon spread over two months, with 8 weekly tactical periods and hourly operational period. Each week is composed of seven days of ten operating hours each. The $\epsilon$-method is implemented with a guaranteed optimality gap of $0.1\%$.

\begin{table}[H]
\centering
\caption{Experimental Instances}
\label{Instances}
\begin{tabular}{ccccc}
Instance & Access hubs & Local cells & Area covered (sq.km) & Population covered  \\
\hline
A        & 39               & 1           & 24.2   & 338,000             \\
B        & 54               & 2           & 42.1   & 590,000             \\
C        & 138              & 4           & 66     & 924,000             \\
D        & 421              & 10          & 178.4  & 2,500,000           \\
E        & 838              & 20          & 410    & 5,740,000
\end{tabular}
\end{table}

As benchmark solutions, static capacity deployments are considered for each instance. Such static capacity deployment represents the minimum capacity module deployment required over the network of access hub locations to satisfy storage requirements for all operational periods within the planning horizon $T$ without being able to update capacity over time or use capacity pooling recourse actions. Benchmark solutions are found by solving $\min  \{ \sum_{t \in T} \sum_{l \in \mathcal{L}\cup \mathcal{W}} h_l S_l(t) \}$ such that $S_l(t) \in \{v S_l(t) \geq D_l(\tau, \omega), \forall l \in L, \tau \in T_{t}, t \in T, \omega \in \Omega \}$ over the entire planning horizon with no relocation or recourse by relaxing spatial constraints to ensure feasibility. \\
An instance has more or less savings potential depending on its demand dynamics and network configuration. Although assessing the potential of capacity pooling a priori is non trivial, the potential of capacity relocation can be assessed by a lower bound to the dynamic capacity deployment problem with no capacity pooling.  Define $\Tilde{S}_l(t)$ as the maximum number of capacity modules required at location $l$ in any operational period associated with tactical period $t$ in all considered scenarios; that is $\Tilde{S}_l(t) = max(\ceil{D_l(\tau, \omega)/v} , \forall \tau \in T_t, \omega \in \Omega)$. Then, an instance's capacity relocation cost savings potential can be computed by factoring in holding costs while ignoring relocation costs, producing a lower bound for the dynamic capacity deployment problem with no capacity pooling, with objective value $\sum_{t \in T} h_l\Tilde{S}_l(t)$. Benchmark solutions and relocation potential for the considered instances are summarized in Table \ref{Benchmark}.

\begin{table}[H]
\centering
\caption{Benchmark Solutions and Relocation Potential}
\label{Benchmark}
\begin{tabular}{cccc}
Instance & Total cost & Capacity & Potential cost savings  \\
\hline
A        &       $\$138,966$         &      189  & $7.93\%$     \\
B        &       $\$181,275$         &      245  & $8.72\%$     \\
C        &       $\$453,137$        &       618  & $7.58\%$    \\
D        &       $\$1,325,490$        &     1820 & $7.29\%$     \\
E        &       $\$2,840,510$        &     3818    & $9.25\%$
\end{tabular}
\end{table}

The initial capacity deployment is defined by running the proposed solution approach for the tactical period immediately preceding the studied planning horizon by relaxing constraint (\ref{IB}).
The default values of input parameters are estimated relying on company experts and presented in \ref{Experiment Parameters}.
Each instance is assigned one depot in one of its local hub locations to store unused capacity modules at no cost. The number of modules available $I_0$ and the penalty cost $p_l$ are set to large values (respectively $5000$ modules and $\$100,000$ per modules in order to prevent full recourse actions by lack of capacity and focus on feasible capacity deployments with capacity pooling. As suggested by \cite{winkenbach2016enabling} (through simulation) when studying a french parcel express company, this paper considers the value of the $k$ constants to be 0.82 for riders and 1.15 for couriers. \\
All experiments were implemented in Python 3.7 using Gurobi 9.0 as the solver and were computed using 40 logical processors on an AMD EPYC Processor @ 2500GHz.

\subsection{Computational Performance}

\begin{figure}[H]
\caption{Benders Algorithm Improvements Comparison on Instance C}
\label{benders performance C}
\centering
\fbox{\includegraphics[scale=0.6]{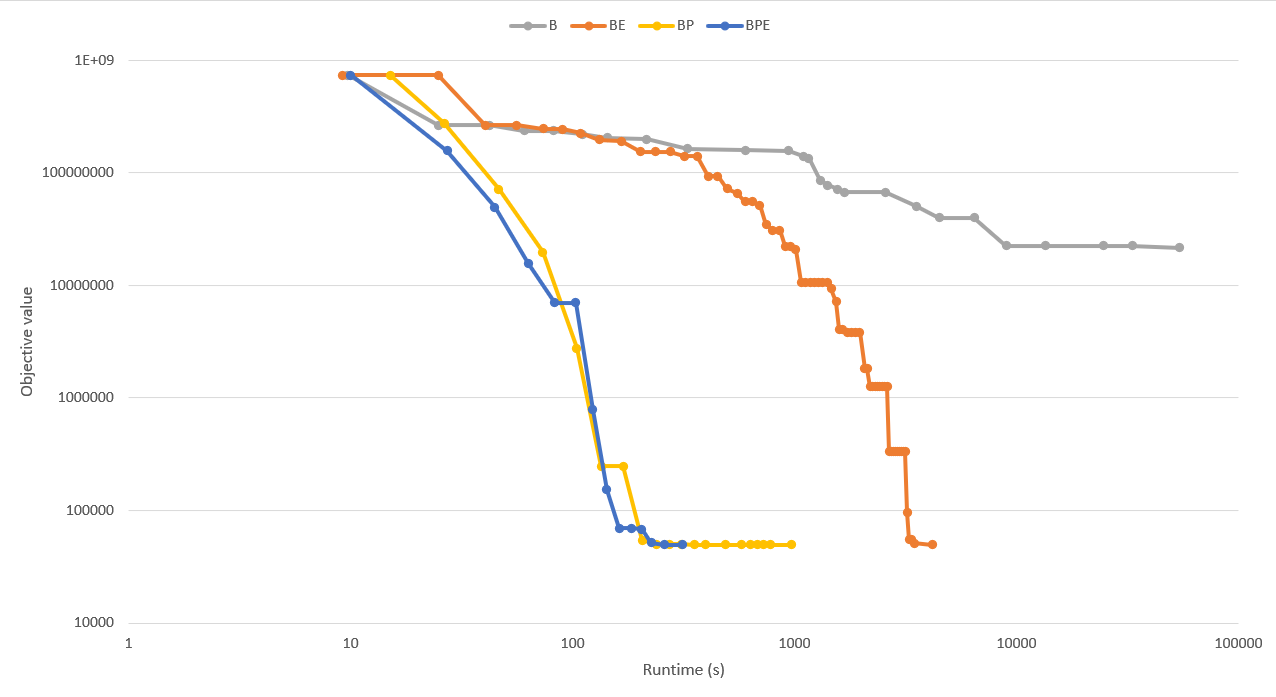}}
\end{figure}

The experiments presented in this section study the computational performance of the proposed solution approach when tackling instances of different sizes. The first experiment aims at validating the efficiency of the proposed acceleration methods in section (\ref{preprocessing and acceleration}) for the Benders algorithm. It examines the impact of combinations of the acceleration methods on the runtime of the Benders algorithm for solving the optimization model (\ref{Obj}-\ref{integrality constraint}) for one relocation period with no lookahead. Figure \ref{benders performance C} display the runtimes for instances C with a capacity pooling distance of 1km and a time cutoff of 15 hours; B represents the original Benders algorithm developed in section (\ref{BD}); BP represents Benders with pareto-optimal cuts; BE represents Benders with the $\epsilon$-optimal method; and BPE represents BP with the $\epsilon$-optimal method.\\ Figure \ref{benders performance C} suggests that pareto-optimal cuts  have the strongest impact on computational performance as it allows the BP algorithm to converge in 965 seconds when the B algorithm did not converge within the time limit. The $\epsilon$-optimal method suggests a significant improvement compared to the original Benders algorithm, and has an advantage over BP when close to optimality (while guaranteeing a solution within $0.1\%$ of optimality). Similar behaviors can be observed for larger instances, with BPE outperforming B, BP and BE.

\begin{figure}[H]
\centering
\fbox{\includegraphics[scale=0.60]{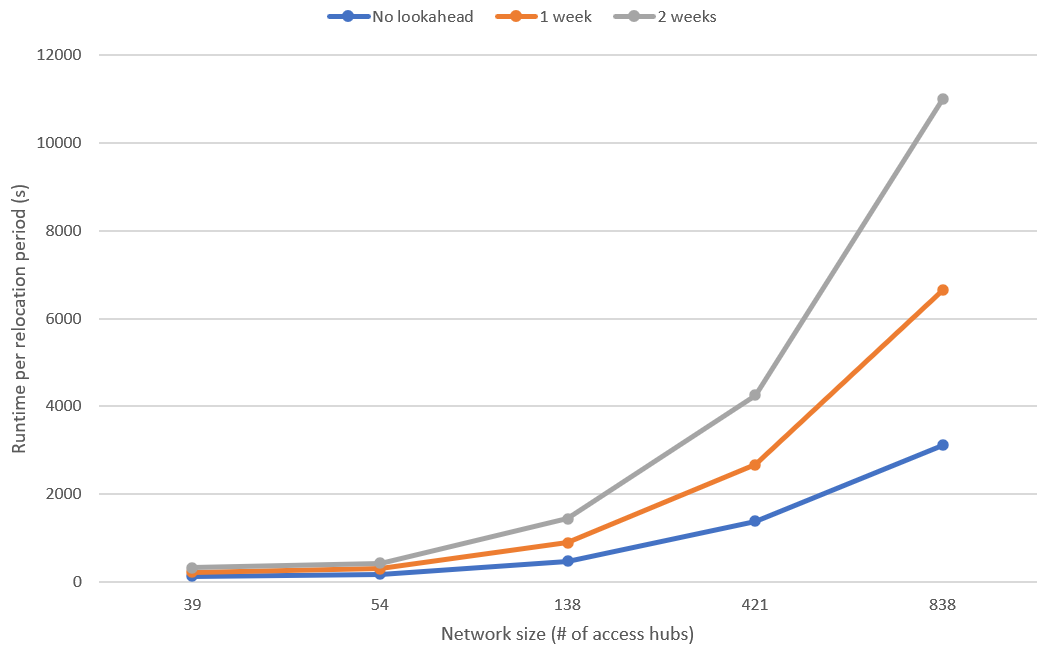}}
\caption{Computational Performance for Different Lookahead Values as a Function of Network Size}
\label{algorithm_perf}
\end{figure}

Next, Figure \ref{algorithm_perf} depicts the computational performance of the proposed solution approach for different lookahead values as a function of network size. Each data point is the average runtime per period for a minimum sample set of 16 instances (8 relocation periods times 2 capacity pooling distances) and a maximum of 48 instances (8 relocation periods times 6 capacity pooling distances) based on the other experiments presented in the paper. The first observation is that the proposed solution approach is efficient in solving large-scale instances considered in this paper (838 access hubs), with a maximum runtime around 3 hours (with 2 weeks lookahead); this result suggests tractability for most urban area sizes, including megacities. The second observation is that adding tactical lookahead reasonably increase runtime: 1 week and 2 weeks lookahead runtimes are respectively at most 2.1 times and 3.5 times as long as no lookahead runtimes wihtin the range of network sizes considered.

\subsection{Comparative Results 
}

The results presented in this section highlight the benefits of relocating capacity dynamically over time and allowing capacity pooling compared to a static capacity deployment with no capacity pooling.
Results are summarized in Table \ref{Result table} for different lookahead values and capacity pooling distance (in km). Table \ref{Result table} presents total costs of the network, deployed capacity (maximum number of modules), relocation share (average number of relocations per period as a share of capacity), and cost and capacity savings with respect to the static counterpart.

First, cost and capacity savings are observed in all the instances.
Maximum cost savings of $28.3\%$ and capacity savings of $26.46\%$ are reached for instance A with a capacity pooling distance of 2km and a 2 weeks tactical lookahead. Most of these savings are a result of the capacity pooling recourse as savings with capacity pooling of 0km indicate a much lower savings (maximum of $6.26\%$ cost savings). Note that for each instance, savings with no capacity pooling are less than potential savings presented in Table \ref{Benchmark} (where relocation costs are not accounted for). The average number of relocations per period represent up to $7.31\%$ of the capacity, and is decreasing as more tactical lookahead is added; capacity deployments are gradually reconfiguring networks. Capacity savings indicate that the total number of modules required (both deployed and stored at a depot) is inferior to the number of modules required in static counterparts. Capacity savings also increase as capacity pooling is available, making the total capital invested in capacity modules inferior than in static counterparts. \\
\begin{table}[H]
\centering
\caption{Results Highlights and Comparison with static capacity deployments}
\label{Result table}
\resizebox{\columnwidth}{!}{
\begin{tabular}{cccccccc}
Instance           & Pooling
  distance & Lookahead ($\theta$)    & Total cost              & Capacity & Relocation share & Cost savings & Capacity savings  \\ 
\hline
\multirow{6}{*}{A} & \multirow{3}{*}{0}      & 0 & \$131,333  & 183    & 7.31\%                         & 5.49\%                           & 3.17\%         \\
                   &                         & 1       & \$130,269  & 187      & 4.88\%                               & 6.26\%                           & 1.06\%                                \\
                   &                         & 2      & \$130,083  & 184                               & 4.82\%                               & 6.39\%                           & 2.65\%                                \\ 
\cline{2-8}
                   & \multirow{3}{*}{2}      & 0 & \$100,925  & 139                 & 7.19\%                              & 27.37\%                          & 26.46\%                               \\
                   &                         & 1       & \$99,779 & 139                  & 6.47\%                               & 28.20\%                          & 26.46\%                               \\
                   &                         & 2      & \$99,644 & 139            & 6.03\%                               & 28.30\%                          & 26.46\%                               \\ 
\hline
\multirow{6}{*}{B} & \multirow{3}{*}{0}      & 0 & \$170,189  & 237          & 6.86\%                              & 6.12\%                           & 3.27\%                                \\
                   &                         & 1       & \$169,940  & 239                & 6.64\%                               & 6.25\%                           & 2.45\%                                \\
                   &                         & 2      & \$169,174  & 239             & 6.33\%                               & 6.68\%                           & 2.45\%                                \\ 
\cline{2-8}
                   & \multirow{3}{*}{2}      & 0 & \$140,470  & 195                    & 6.60\%                               & 22.51\%                          & 20.41\%                               \\
                   &                         & 1       & \$138,847  & 195          & 6.15\%                               & 23.41\%                          & 20.41\%                               \\
                   &                         & 2      & \$138,841  & 195                    & 6.15\%                               & 23.41\%                          & 20.41\%                               \\ 
\hline
\multirow{6}{*}{C} & \multirow{3}{*}{0}      & 0 & \$429,967  & 594         & 6.29\%                               & 5.11\%                           & 3.88\%                                \\
                   &                         & 1       & \$428,683  & 595                  & 6.11\%                               & 5.40\%                           & 3.72\%                                \\
                   &                         & 2      & \$428,544  & 595           & 5.95\%                               & 5.43\%                           & 3.72\%                                \\ 
\cline{2-8}
                   & \multirow{3}{*}{2}      & 0 & \$400,482  & 548              & 6.20\%                               & 11.62\%                          & 11.33\%                               \\
                   &                         & 1       & \$396,944  & 550             & 5.64\%                               & 12.40\%                          & 11.00\%                               \\
                   &                         & 2      & \$396,776  & 550           & 5.59\%                               & 12.44\%                          & 11.00\%                               \\ 
\hline
\multirow{6}{*}{D} & \multirow{3}{*}{0}      & 0 & \$1,263,680 & 1746         & 7.01\%                               & 4.66\%                           & 4.07\%                                \\
                   &                         & 1       & \$1,255,070 & 1747           & 6.18\%                               & 5.31\%                           & 4.01\%                                \\
                   &                         & 2      & \$1,249,130 & 1737       & 6.02\%                               & 5.76\%                           & 4.56\%                                \\ 
\cline{2-8}
                   & \multirow{3}{*}{2}      & 0 & \$1,228,360 & 1698           & 6.39\%                               & 7.33\%                           & 6.70\%                                \\
                   &                         & 1       & \$1,219,230 & 1699           & 6.02\%                               & 8.02\%                           & 6.65\%                                \\
                   &                         & 2      & \$1,219,050 & 1696           & 6.01\%                               & 8.03\%                           & 6.81\%                                \\ 
\hline
\multirow{6}{*}{E} & \multirow{3}{*}{0}      & 0 & \$2,646,090 & 3632          & 7.00\%                               & 6.84\%                           & 4.87\%                                \\
                   &                         & 1       & \$2,624,860 & 3626             & 6.52\%                               & 7.59\%                           & 5.03\%                                \\
                   &                         & 2      & \$2,624,380 & 3626       & 6.51\%                              & 7.61\%                           & 5.03\%                                \\ 
\cline{2-8}
                   & \multirow{3}{*}{2}      & 0 & \$2,614,330 & 3587       & 6.93\%                               & 7.96\%                           & 6.05\%                                \\
                   &                         & 1       & \$2,596,690 & 3587         & 6.55\%                               & 8.58\%                           & 6.05\%                                \\
                   &                         & 2      & \$2,595,550 & 3589          & 6.53\%                               & 8.62\%                           & 6.00\%                               
\end{tabular}
}
\end{table}

Furthermore, the results show that adding tactical lookahead is beneficial for all instances with and without capacity pooling by improving cost savings and decreasing the number of relocations. The role of tactical lookahead is to anticipate future needs and avoid relocations that will be reverted to in the future. Lookahead can be seen as the flexibility hedging of the solution approach to avoid relocations under uncertainty. However, the difference between one week and two weeks of tactical lookahead is more subtle with smaller cost improvements. These results suggest that solution's quality increase with lookahead ($\theta$), offering extra cost savings. Tactical lookahead anticipates for future relocations therefore decreasing relocation share at the cost of slightly higher capacity deployments. However, there does not seem to be significant improvements from extending the lookahead from one week to two weeks, especially when considering the additional computational runtime.\\
Lastly, capacity pooling brings significant value to instance A, B, and C, but less cost savings improvements for instance D and E. This is probably due to the the fact that instances D and E have lower hub density, increasing the distance between access hubs (see Table \ref{Instances}). Section \ref{pooling varation} examine the impact of capacity pooling distance in more details by focusing on instance C.
\subsection{Capacity Pooling Variations} \label{pooling varation}

\begin{figure}[H]
\centering
\fbox{\includegraphics[scale=1]{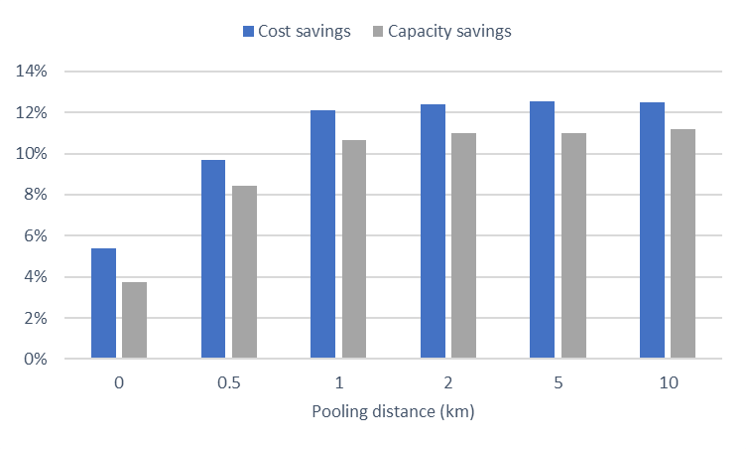}}
\caption{Cost and Capacity Savings as a Function of Pooling Distance}
\label{fig: Pooling Exp}
\end{figure}

This experiment examines the effect of capacity pooling as a way to further decrease costs. Table \ref{table: capacity pooling} summarizes the effect of different capacity pooling distances (in km) on instance C's solutions. It presents average additional rider and courier travel (induced by detours), and cost and capacity savings for instance C. Figure \ref{fig: Pooling Exp} displays a plot of cost and capacity savings as a function of pooling distance.

\begin{table} [h]
\centering
\caption{Sensitivity Analysis on Capacity Pooling Distance for Instance C ($\theta = 1$)}
\label{table: capacity pooling}
\begin{tabular}{ccccc}
Pooling distance  & Rider travel & Courier travel & Cost savings & Capacity savings  \\ 
\hline
0               & 0.00                    & 0.00                      & 5.40\%       & 3.72\%           \\ 

0.5             & 0.35                    & 4.59                      & 9.67\%      & 8.41\%           \\
                    
1             & 0.49                    & 9.86                      & 12.12\%      & 10.68\%           \\
                    
2             & 0.48                    & 9.71                      & 12.40\%      & 11.00\%           \\
                    
5             & 0.48                    & 9.71                      & 12.53\%      & 11.00\%           \\
                    
10            & 0.48                    & 9.71                      & 12.51\%      & 11.17\%           \\
\end{tabular}
\end{table}

The increase in capacity pooling distance allows to produce superior solutions but only until a maximum of $12.55\%$ is reached with a pooling distance of 5km. This trend can clearly be seen in Figure \ref{fig: Pooling Exp}. Indeed, no matter how large capacity pooling pooling neighborhoods are, constraints (\ref{Rider_Synchro_Constraint}) and (\ref{Courier_Synchro_Constraint}) limit capacity poolings from an operational point of view: riders and couriers cannot perform long distance detours as it would disrupt their activity by delaying other pickup / deliveries. Table \ref{pooling varation} shows that most of the additional travel induced by detours is performed by couriers; since riders have larger carrying capacity, one rider detour may require multiple courier detours. Note also that since couriers are often using lightweight vehicles (if any vehicles at all), long distance detours may not be practical which may also limit the capacity pooling distance from a design perspective. The same behavior can be observed for the other instances.

\subsection{Holding Costs Versus Relocation Costs}

This experiment examines the influence of relocation costs and holding costs on dynamic capacity deployments. Intuitively, two extreme cases can be identified: (1) if holding costs are negligible compared to relocation costs, there is no incentive to dynamically adjust capacity, and (2) if relocation costs are negligible compared to holding costs, a myopic view of the problem would be optimal as anticipating future relocations does not save cost. Apart from extreme cases, variations of holding costs and relocation cost can represent different urban environments. A very dense city may have high holding costs (prime real estate) and low relocation costs (short distances between locations). In this experiment, four cases are examined: High-high, High-low, Low-high and Low-low, where High and high respectively represent high holding costs and high relocation costs and Low and low respectively represent low holding costs and low relocation costs. Cost vectors are scaled linearly and high costs are chosen to be $100\%$ higher than baseline values while low costs are assumed to be $50\%$ lower than baseline values. Table \ref{table: H vs C} presents total cost, capacity (maximum number of modules deployed), relocations, relocation share (average number of relocations per period as a share of capacity), and cost and capacity savings for instance C with a capacity pooling distance of 2km. Savings are computed comparing to benchmark solutions with corresponding cost adjustements (holding costs).

\begin{table} [H]
\centering
\caption{Impact of Holding versus Relocation Costs on Instance C}
\label{table: H vs C}
\begin{tabular}{c c c c}
Case      & Relocation share & Cost savings & Capacity savings  \\
\hline
High-low  & $6.12\%$                      & $13.07\%$      & $9.62\%$            \\

High-high & $5.22\%$                      & $11.76\%$      & $9.78\%$            \\

Low-low   & $6.20\%$                      & $13.10\%$      & $8.13\%$           \\

Low-high  & $3.73\%$                      & $9.02\%$       & $7.79\%$           \\

\end{tabular}
\end{table}





A first observation is that cases where relocation costs are low perform best with costs savings around $13.1\%$, regardless of holding costs. When relocation costs are high, cost savings are worse, especially when holding costs are low ($9.02\%$).
Low holding cost cases deploy more capacity modules overall which impact capacity savings, but are still able to reach high cost savings when relocation costs are low. The combination of low holding costs and high relocation costs decreases opportunities for worthy relocations (only $3.73\%$ relocation share), requiring more capacity deployed and therefore limiting cost savings ($9.02\%$). Similar behavior can be observed on the other instances. \\
Overall, this experiment indicates that denser urban environment (high holding costs) tend to be better candidates for dynamic capacity management of access hub networks. Moreover, low relocation costs (i.e. easy installation and good mobility of capacity modules) can make any urban environment a worthy candidate for such capacity management strategy. Finally, the combination of lesser dense urban environment and high relocation costs significantly limits opportunity for cost savings.

\section{Conclusion} \label{Conclusion}

This paper defines and formulates the Dynamic Pooled Capacity Deployment Problem in the context of urban parcel logistics. This problem involves a tactical decision on the relocation of capacity modules over a network of discrete locations associated with stochastic demand requirements. To improve the quality of the capacity deployment decisions, the proposed model integrates an estimate of the difference of operations cost, which includes capacity assignment decisions with the possibility of capacity pooling between neighboring locations. It also integrates synchronization requirements of the 2-echelon routing subproblems, using an analytical derivate from the route length estimation function. The dynamic problem is modeled and approximated with a two-stage stochastic program with recourse, where all capacity deployment decisions on a finite planning horizon are moved to the first stage. Due to the uncertainty of capacity requirements and the challenges of solving the MIP formulation for realistic networks of several hundreds locations, a roll-out approach with lookahead based on a Benders decomposition of the finite planning horizon problem coupled with acceleration methods is proposed. Five instances of networks of different sizes are presented to perform computational experiments to test the performance of the proposed approach and assess the potential of the defined capacity deployment strategy. \\
Results show that the proposed approach produces solutions in a reasonable time even for large scale instances of up to 838 hubs. They suggest that a dynamic capacity deployment strategy with capacity pooling has a significant advantage over a static capacity deployment strategy for access hub networks, with up to $28\%$ cost savings and $26\%$ capacity savings. Results also show that one-week lookahead helps producing superior solutions by anticipating future relocations, but adding a two-weeks lookahead does not make a significant improvement. Increasing the capacity pooling distance, while increasing computing time, tend to increase opportunities for cost savings by allowing more locations to pool capacity until an operational feasibility threshold is reached. Dynamically adjusting workforce assignment in the network was not explored but could potentially overcome this limitation. Denser urban environments (i.e. with higher real estate costs) are natural candidates for dynamic capacity deployments as relocation costs are more easily overcome by holding costs. However, relocation costs are the most limiting when it comes to cost savings. Technology solutions featuring cheaper installation costs and high degree of mobility make it more interesting to consider periodic network reconfigurations. \\
The implementation of such innovation also has management challenges not studied in this paper. For instance, implementation may require a more agile workforce, specialized training and targeted hiring enabling a data-driven approach to managing network capacity. Management challenges also need to be considered by decision makers along with the potential reduction of fixed-assets offered by capacity savings when evaluating the solution for implementation.\\
Finally, there are numerous research avenues around reconfigurable networks, dynamic capacity management and access hubs in urban parcel logistics. Where technology allows for very frequent network reconfiguration, solutions featuring not only modular but mobile capacity (e.g. on wheels) and near real-time capacity relocation can become relevant as a complement to the proposed dynamic capacity deployment strategy. Moreover, the possibility of updating operations planning as needed (e.g. dynamic routing, dynamic staffing) can unlock the potential of capacity pooling not only as a recourse but as an integral part of network design and operations planning.

\bibliographystyle{informs2014trsc}

\bibliography{Bibliography}

\begin{thebibliography}{61}
\providecommand{\natexlab}[1]{#1}
\providecommand{\url}[1]{\texttt{#1}}
\providecommand{\urlprefix}{URL }

\bibitem[{Aghezzaf(2005)}]{aghezzaf2005capacity}
Aghezzaf E, 2005 \emph{Capacity planning and warehouse location in supply
  chains with uncertain demands}. \emph{Journal of the Operational Research
  Society} 56(4):453--462.

\bibitem[{Anderluh, Hemmelmayr, \protect\BIBand{}
  Nolz(2017)}]{anderluh2017synchronizing}
Anderluh A, Hemmelmayr V, Nolz P, 2017 \emph{Synchronizing vans and cargo bikes
  in a city distribution network}. \emph{Central European Journal of Operations
  Research} 25(2):345--376.

\bibitem[{Ansari et~al.(2018)Ansari, Ba{\c{s}}dere, Li, Ouyang,
  \protect\BIBand{} Smilowitz}]{ansari2018advancements}
Ansari S, Ba{\c{s}}dere M, Li X, Ouyang Y, Smilowitz K, 2018 \emph{Advancements
  in continuous approximation models for logistics and transportation systems:
  1996--2016}. \emph{Transportation Research Part B: Methodological}
  107:229--252.

\bibitem[{Arabani \protect\BIBand{} Farahani(2012)}]{arabani2012facility}
Arabani A, Farahani R, 2012 \emph{Facility location dynamics: An overview of
  classifications and applications}. \emph{Computers \& Industrial Engineering}
  62(1):408--420.

\bibitem[{Arvidsson \protect\BIBand{} Pazirandeh(2017)}]{arvidsson2017ex}
Arvidsson N, Pazirandeh A, 2017 \emph{An ex ante evaluation of mobile depots in
  cities: A sustainability perspective}. \emph{International Journal of
  Sustainable Transportation} 11(8):623--632.

\bibitem[{Benders(2005)}]{Benders2005}
Benders J, 2005 \emph{Partitioning procedures for solving mixed-variables
  programming problems}. \emph{Computational Management Science} 2(1):3--19.

\bibitem[{Benjelloun \protect\BIBand{} Crainic(2008)}]{BenjellounTrends}
Benjelloun A, Crainic T, 2008 \emph{Trends, challenges, and perspectives in
  city logistics}. \emph{Transportation and land use interaction, proceedings
  TRANSLU} 8:269--284.

\bibitem[{Bergmann, Wagner, \protect\BIBand{}
  Winkenbach(2020)}]{bergmann2020integrating}
Bergmann FM, Wagner SM, Winkenbach M, 2020 \emph{Integrating first-mile pickup
  and last-mile delivery on shared vehicle routes for efficient urban
  e-commerce distribution}. \emph{Transportation Research Part B:
  Methodological} 131:26--62.

\bibitem[{{CityLog}(2012)}]{CityLogBentoBox}
{CityLog}, 2012 \emph{Citylog project bentobox}.
  \urlprefix\url{https://www.berlin.de/senuvk/verkehr/politik_planung/projekte/citylog/index.shtml},
  [Online; accessed March 10, 2020].

\bibitem[{Crainic \protect\BIBand{} Montreuil(2016)}]{crainic2016physical}
Crainic T, Montreuil B, 2016 \emph{Physical internet enabled hyperconnected
  city logistics}. \emph{Transportation Research Procedia} 12:383--398.

\bibitem[{Crainic et~al.(2016)Crainic, Hewitt, Toulouse, \protect\BIBand{}
  Vu}]{crainic2016service}
Crainic TG, Hewitt M, Toulouse M, Vu DM, 2016 \emph{Service network design with
  resource constraints}. \emph{Transportation Science} 50(4):1380--1393.

\bibitem[{Crainic, Ricciardi, \protect\BIBand{}
  Storchi(2004)}]{crainic2004advanced}
Crainic TG, Ricciardi N, Storchi G, 2004 \emph{Advanced freight transportation
  systems for congested urban areas}. \emph{Transportation Research Part C:
  Emerging Technologies} 12(2):119--137.

\bibitem[{Cuda, Guastaroba, \protect\BIBand{} Speranza(2015)}]{Cuda2015}
Cuda R, Guastaroba G, Speranza MG, 2015 \emph{A survey on two-echelon routing
  problems}. \emph{Computers \& Operations Research} 55:185--199.

\bibitem[{Daganzo(1984)}]{Daganzo94}
Daganzo C, 1984 \emph{The length of tours in zones of different shapes}.
  \emph{Transportation Research Part B: Methodological} 18(2):135--145.

\bibitem[{Daganzo(2005)}]{DaganzoBook}
Daganzo C, 2005 \emph{Logistics systems analysis} (Springer Science \& Business
  Media), ISBN 3540275169.

\bibitem[{Drexl \protect\BIBand{} Schneider(2015)}]{drexl2015survey}
Drexl M, Schneider M, 2015 \emph{A survey of variants and extensions of the
  location-routing problem}. \emph{European Journal of Operational Research}
  241(2):283--308.

\bibitem[{eMarketer(2018)}]{statista_2019_speed}
eMarketer, 2018 \emph{Items that internet users in the united states want to
  have vs. have received via same-day delivery as of february 2018}.
  \urlprefix\url{https://www.statista.com/statistics/861389/items-internet-users-want-have-vs-received-same-day-delivery/}.

\bibitem[{eMarketer(2019)}]{statista_2019_growth}
eMarketer, 2019 \emph{Annual retail e-commerce sales growth worldwide from 2014
  to 2023}.
  \urlprefix\url{https://www.statista.com/statistics/288487/forecast-of-global-b2c-e-commerce-growt/}.

\bibitem[{Enthoven et~al.(2020)Enthoven, Jargalsaikhan, Roodbergen, uit~het
  Broek, \protect\BIBand{} Schrotenboer}]{enthoven2020two}
Enthoven D, Jargalsaikhan B, Roodbergen K, uit~het Broek M, Schrotenboer A,
  2020 \emph{The two-echelon vehicle routing problem with covering options:
  City logistics with cargo bikes and parcel lockers}. \emph{Computers \&
  Operations Research} 118:104919.

\bibitem[{Erera(2000)}]{EreraThesis}
Erera A, 2000 \emph{Design of large-scale logistics systems for uncertain
  environments} (University of California, Berkeley).

\bibitem[{Faug\`ere \protect\BIBand{} Montreuil(2020)}]{FaugereCIE}
Faug\`ere L, Montreuil B, 2020 \emph{Smart locker bank design optimization for
  urban omnichannel logistics: Assessing monolithic vs. modular
  configurations}. \emph{Computers \& Industrial Engineering} 139:105544.

\bibitem[{Franceschetti et~al.(2017)Franceschetti, Honhon, Laporte,
  Van~Woensel, \protect\BIBand{} Fransoo}]{franceschetti2017strategic}
Franceschetti A, Honhon D, Laporte G, Van~Woensel T, Fransoo JC, 2017
  \emph{Strategic fleet planning for city logistics}. \emph{Transportation
  Research Part B: Methodological} 95:19--40.

\bibitem[{Francis, Smilowitz, \protect\BIBand{} Tzur(2008)}]{francis2008period}
Francis P, Smilowitz K, Tzur M, 2008 \emph{The period vehicle routing problem
  and its extensions}. \emph{The vehicle routing problem: latest advances and
  new challenges}, 73--102 (Springer).

\bibitem[{Fransoo, Blanco, \protect\BIBand{}
  Argueta(2017)}]{fransoo2017reaching}
Fransoo J, Blanco E, Argueta CM, 2017 \emph{Reaching 50 million nanostores:
  retail distribution in emerging megacities} (CreateSpace Independent
  Publishing Platform).

\bibitem[{Ghiani, Guerriero, \protect\BIBand{}
  Musmanno(2002)}]{ghiani2002capacitated}
Ghiani G, Guerriero F, Musmanno R, 2002 \emph{The capacitated plant location
  problem with multiple facilities in the same site}. \emph{Computers \&
  Operations Research} 29(13):1903--1912.

\bibitem[{Gonzalez-Feliu(2009)}]{gonzalez2009n}
Gonzalez-Feliu J, 2009 \emph{The n-echelon location routing problem: concepts
  and methods for tactical and operational planning}. \emph{Technical Report
  halshs-00422492, Laboratoire d ‘Economie des Transports, Institut des
  Sciences de l‘Homme} .

\bibitem[{Hewitt et~al.(2019)Hewitt, Crainic, Nowak, \protect\BIBand{}
  Rei}]{hewitt2019scheduled}
Hewitt M, Crainic TG, Nowak M, Rei W, 2019 \emph{Scheduled service network
  design with resource acquisition and management under uncertainty}.
  \emph{Transportation Research Part B: Methodological} 128:324--343.

\bibitem[{Janjevic \protect\BIBand{} Ndiaye(2014)}]{janjevic_development_2014}
Janjevic M, Ndiaye AB, 2014 \emph{Development and {Application} of a
  {Transferability} {Framework} for {Micro}-consolidation {Schemes} in {Urban}
  {Freight} {Transport}}. \emph{Procedia - Social and Behavioral Sciences}
  125:284--296,
  \urlprefix\url{http://dx.doi.org/10.1016/j.sbspro.2014.01.1474}.

\bibitem[{Janjevic \protect\BIBand{}
  Winkenbach(2020)}]{janjevic2020characterizing}
Janjevic M, Winkenbach M, 2020 \emph{Characterizing urban last-mile
  distribution strategies in mature and emerging e-commerce markets}.
  \emph{Transportation Research Part A: Policy and Practice} 133:164--196.

\bibitem[{Janjevic, Winkenbach, \protect\BIBand{}
  Merch{\'a}n(2019)}]{janjevic2019integrating}
Janjevic M, Winkenbach M, Merch{\'a}n D, 2019 \emph{Integrating
  collection-and-delivery points in the strategic design of urban last-mile
  e-commerce distribution networks}. \emph{Transportation Research Part E:
  Logistics and Transportation Review} 131:37--67.

\bibitem[{Jena, Cordeau, \protect\BIBand{} Gendron(2015)}]{jena2015dynamic}
Jena S, Cordeau J, Gendron B, 2015 \emph{Dynamic facility location with
  generalized modular capacities}. \emph{Transportation Science}
  49(3):484--499.

\bibitem[{Klibi, Martel, \protect\BIBand{} Guitouni(2016)}]{Klibi2016}
Klibi W, Martel A, Guitouni A, 2016 \emph{The impact of operations
  anticipations on the quality of stochastic location-allocation models}.
  \emph{Omega} 62:19--33.

\bibitem[{Leonardi, Browne, \protect\BIBand{} Allen(2012)}]{leonardi2012before}
Leonardi J, Browne M, Allen J, 2012 \emph{Before-after assessment of a
  logistics trial with clean urban freight vehicles: A case study in london}.
  \emph{Procedia-Social and Behavioral Sciences} 39:146--157.

\bibitem[{Maheo, Kilby, \protect\BIBand{}
  Van~Hentenryck(2019)}]{maheo2019benders}
Maheo A, Kilby P, Van~Hentenryck P, 2019 \emph{Benders decomposition for the
  design of a hub and shuttle public transit system}. \emph{Transportation
  Science} 53(1):77--88.

\bibitem[{Malladi, Erera, \protect\BIBand{}
  White~III(2020)}]{malladiererawhite}
Malladi SS, Erera AL, White~III CC, 2020 \emph{A dynamic mobile production
  capacity and inventory control problem}. \emph{Iise Transactions}
  52(8):926--943.

\bibitem[{Mancini(2013)}]{mancini2013multi}
Mancini S, 2013 \emph{Multi-echelon distribution systems in city logistics}.
  \emph{European Transport} 54 (2):1--24.

\bibitem[{Manzini \protect\BIBand{} Gebennini(2008)}]{manzini2008optimization}
Manzini R, Gebennini E, 2008 \emph{Optimization models for the dynamic facility
  location and allocation problem}. \emph{International Journal of Production
  Research} 46(8):2061--2086.

\bibitem[{Marcotte \protect\BIBand{} Montreuil(2016)}]{marcotte2016introducing}
Marcotte S, Montreuil B, 2016 \emph{Introducing the concept of hyperconnected
  mobile production}. \emph{Progress in Material Handling Research} .

\bibitem[{Marcotte, Montreuil, \protect\BIBand{}
  Coelho(2015)}]{marcotte2015modeling}
Marcotte S, Montreuil B, Coelho L, 2015 \emph{Modeling of physical internet
  enabled interconnected modular production}. \emph{Proceedings of 2nd
  International Physical Internet Conference, Paris, France}.

\bibitem[{Marujo et~al.(2018)Marujo, Goes, D'Agosto, Ferreira, Winkenbach,
  \protect\BIBand{} Bandeira}]{marujo2018assessing}
Marujo L, Goes G, D'Agosto M, Ferreira A, Winkenbach M, Bandeira R, 2018
  \emph{Assessing the sustainability of mobile depots: The case of urban
  freight distribution in rio de janeiro}. \emph{Transportation Research Part
  D: Transport and Environment} 62:256--267.

\bibitem[{Melo, Nickel, \protect\BIBand{} Da~Gama(2006)}]{melo2006dynamic}
Melo M, Nickel S, Da~Gama F, 2006 \emph{Dynamic multi-commodity capacitated
  facility location: a mathematical modeling framework for strategic supply
  chain planning}. \emph{Computers \& Operations Research} 33(1):181--208.

\bibitem[{Montreuil(2011)}]{montreuil2011toward}
Montreuil B, 2011 \emph{Toward a physical internet: meeting the global
  logistics sustainability grand challenge}. \emph{Logistics Research}
  3(2-3):71--87.

\bibitem[{Montreuil et~al.(2018)Montreuil, Buckley, Faugere, Khir,
  \protect\BIBand{} Derhami}]{montreuil2018IMHRC}
Montreuil B, Buckley S, Faugere L, Khir R, Derhami S, 2018 \emph{Parcel
  logistics hub design: The impact of modularity and hyperconnectivity}.
  \emph{Progress in Material Handling Research} .

\bibitem[{Papadakos(2009)}]{papadakos2009integrated}
Papadakos N, 2009 \emph{Integrated airline scheduling}. \emph{Computers \&
  Operations Research} 36(1):176--195.

\bibitem[{Perboli, Tadei, \protect\BIBand{} Vigo(2011)}]{perboli2011two}
Perboli G, Tadei R, Vigo D, 2011 \emph{The two-echelon capacitated vehicle
  routing problem: models and math-based heuristics}. \emph{Transportation
  Science} 45(3):364--380.

\bibitem[{Powell(2007)}]{powell2007approximate}
Powell WB, 2007 \emph{Approximate Dynamic Programming: Solving the curses of
  dimensionality}, volume 703 (John Wiley \& Sons).

\bibitem[{Powell(2019)}]{powell2019unified}
Powell WB, 2019 \emph{A unified framework for stochastic optimization}.
  \emph{European Journal of Operational Research} 275(3):795--821.

\bibitem[{Rahmaniani et~al.(2017)Rahmaniani, Crainic, Gendreau,
  \protect\BIBand{} Rei}]{rahmaniani2017benders}
Rahmaniani R, Crainic T, Gendreau M, Rei W, 2017 \emph{The benders
  decomposition algorithm: A literature review}. \emph{European Journal of
  Operational Research} 259(3):801--817.

\bibitem[{Rais, Alvelos, \protect\BIBand{}
  Carvalho(2014)}]{RaisPandDwithTransshipment}
Rais A, Alvelos F, Carvalho M, 2014 \emph{New mixed integer-programming model
  for the pickup-and-delivery problem with transshipment}. \emph{European
  Journal of Operational Research} 235(3):530--539.

\bibitem[{Savelsbergh \protect\BIBand{}
  Van~Woensel(2016)}]{savelsbergh201650th}
Savelsbergh M, Van~Woensel T, 2016 \emph{50th anniversary invited
  article—city logistics: Challenges and opportunities}. \emph{Transportation
  Science} 50(2):579--590.

\bibitem[{Sch\"ultz(2003)}]{Schultz2003}
Sch\"ultz R, 2003 \emph{Stochastic programming with integer variables}.
  \emph{Mathematical Programming} 97(1):285--309.

\bibitem[{Shapiro, Dentcheva, \protect\BIBand{}
  Ruszczynski(2009)}]{Shapiro2009}
Shapiro A, Dentcheva D, Ruszczynski A, 2009 \emph{Lectures on stochastic
  programming: modeling and theory} (The Society for Industrial and Applied
  Mathematics and the Mathematical Programming Society, Philadelphia, USA).

\bibitem[{Smilowitz \protect\BIBand{} Daganzo(2007)}]{SmilowitzDaganzo}
Smilowitz K, Daganzo C, 2007 \emph{Continuum approximation techniques for the
  design of integrated package distribution systems}. \emph{Networks: An
  International Journal} 50(3):183--196.

\bibitem[{Stodick \protect\BIBand{} Deckert(2019)}]{stodick2019sustainable}
Stodick K, Deckert C, 2019 \emph{Sustainable parcel delivery in urban areas
  with micro de-pots}. \emph{Mobility in a Globalised World 2018} 22:233.

\bibitem[{{The Hub Company}(2019)}]{HubCompany}
{The Hub Company}, 2019 \emph{Dutch hub company presents micro hub solutions
  for city logistics}.
  \urlprefix\url{http://www.citylogistics.info/business/dutch-hub-company-presents-micro-hub-solutions-for-city-logistics/},
  [Online; accessed March 10, 2020].

\bibitem[{{TNT Express}(2013)}]{TNT:MobileAH}
{TNT Express}, 2013 \emph{Tnt express is using a mobile depot in brussels to
  enable more efficient delivery}.
  \urlprefix\url{http://www.tnt.de/__C1257442002D0760.nsf/html/pressemitteilungen_tntexpresseroeffnetmobilesdepotinbruessel.html},
  [Online; accessed March 10, 2020].

\bibitem[{{United Nations}(2018)}]{united20182018}
{United Nations}, 2018 \emph{2018 revision of world urbanization prospects}.

\bibitem[{{UPS}(2018)}]{UPSMicroHub}
{UPS}, 2018 \emph{Letzte meile: Studie zu micro-hubs}.
  \urlprefix\url{https://logistik-heute.de/news/letzte-meile-studie-zu-micro-hubs-14471.html},
  [Online; accessed March 10, 2020].

\bibitem[{Verlinde et~al.(2014)Verlinde, Macharis, Milan, \protect\BIBand{}
  Kin}]{verlinde2014does}
Verlinde S, Macharis C, Milan L, Kin B, 2014 \emph{Does a mobile depot make
  urban deliveries faster, more sustainable and more economically viable:
  results of a pilot test in brussels}. \emph{Transportation Research Procedia}
  4:361--373.

\bibitem[{Wade \protect\BIBand{} Bjerkan(2020)}]{mitsloancovid19}
Wade M, Bjerkan H, 2020 \emph{Three proactive response strategies to covid 19
  business challenges}. \emph{MIT Sloan management review} .

\bibitem[{Winkenbach, Kleindorfer, \protect\BIBand{}
  Spinler(2016)}]{winkenbach2016enabling}
Winkenbach M, Kleindorfer PR, Spinler S, 2016 \emph{Enabling urban logistics
  services at la poste through multi-echelon location-routing}.
  \emph{Transportation Science} 50(2):520--540.

\end{thebibliography}

\appendix

\section{Calibration and Data} \label{Experiment Parameters}

\subsection{In-sample Variability}

In-sample variability was tested with no lookahead for instance A with capacity pooling limited to 1 km for 10 samples. Results are presented in Table \ref{Statistical Analysis}. Coefficient of variation represent the ratio between the standard deviation and the average of solutions' total cost. Statistical gap represent the ratio $(UB-LB)/LB$ where UB and LB are respectively the highest and lowest total cost in the sample.

\begin{table}[h]
\centering
\caption{In-sample Statistical Analysis}
\label{Statistical Analysis}
\begin{tabular}{|c|c|c|c|c|c|c|c|c|}
\hline
Number of scenarios      & 5       & 10     & 20     & 30     & 50     & 75     & 100  & 200     \\
\hline
Coefficient of variation & $3.47\%$  & $2.87\%$ & $2.70\%$ & $2.11\%$ & $1.34\%$ & $1.04\%$ & $0.52\%$ & $0.41\%$  \\
\hline
Statistical gap          & $10.04\%$ & $9.16\%$ & $8.57\%$ & $7.85\%$ & $4.25\%$ & $3.51\%$ & $1.53\%$  & $1.22\%$  \\
\hline
\end{tabular}
\end{table}

\subsection{Cost Estimates}

The capacity module relocation costs $r_a$ include an operational cost of $\$1.50$ per kilometer, and a fixed cost of two operators for two hours at a rate of $\$10$ per hour to uninstall/install modules once at the desired locations:
\begin{subequations}
\begin{align*}
    r_a = 1.50 d_a + 40, \forall a = (i,j) \in A
\end{align*}
\end{subequations}
\noindent Where $d_a$ is the distance between location $i$ and $j$ such that $a = (i,j)$. \\
The holding costs are computed from an amortized acquisition cost of $\$2000$ over 5 years (52 weeks long years), and from a rent cost of $\$75$ per square meter times a location specific factor $(1+f_l)$ randomly generated to represent the real estate difference between locations.
\begin{subequations}
\begin{align*}
    h_l = \frac{2000}{5*52} + 75 (1+ f_l), \forall l \in L
\end{align*}
\end{subequations}
\noindent Where $f_l$ is randomly generated from a uniform distribution over $[2\%, 15\%]$. It is also assumed that modules do not depreciate when stored at depots ($h_l = 0, \forall l \in D$). \\

\subsection{Other Input Parameters}

\begin{table} [h]
\centering
\begin{tabular}{cc|cc}
Parameter      & Value                    & Parameter                   & Value       \\
\hline
$c^C_v$      & $\$1/km$   & $s^R_0$                   & $50 \textit{ km/h}$     \\
$c^C_{v0}$ & $\$0.8/km$ & $\hat{S_l}$ & $15 \textit{ modules}$     \\
$c^R_f$      & $\$10$     &  $t^C_a$                   & $1 \textit{ min}$    \\
$c^R_v$      & $\$1.8/km$ & $t^C_u$                   & $2 \textit{ min}$  \\
$c^R_{v0}$ & $\$1.2/km$ & $t^R_a$                   & $5 \textit{ min}$       \\
$k^C$         & $1.15$                     & $t^R_u$                   & $1 \textit{ min}$       \\
$k^R$         & $0.82$                     & $t^R_s$                   & $5 \textit{ min}$       \\
$p_l$           & $100000/ module$            & $v$                   & $0.75 \textit{ } m^2$       \\
$s^C$         & $7 \textit{ km/h}$                   & $v^C$                   & $0.48 \textit{ } m^2$      \\
$s^C_0$                   & $15 \textit{ km/h}$ &  $v^R$                   & $6.40 \textit{ } m^2$\\
$s^R$                      & $30 \textit{ km/h}$ & &
\end{tabular}
\end{table}


\end{document}